\documentclass[preprint,12pt]{elsarticle}

\usepackage{slashbox} 
\usepackage{mathrsfs}
\usepackage{subfig}
\usepackage{amssymb,amsmath}%
\usepackage{gensymb}
\usepackage{amsthm,bm}%
\usepackage[width=17cm,height=23.5cm,centering]{geometry}%
\usepackage{color}
\usepackage{ulem}
\usepackage{stmaryrd} 
\usepackage{cases,empheq,algorithm,enumitem}
\usepackage[colorlinks,citecolor= blue,linkcolor= blue]{hyperref}
\usepackage{todonotes}  

\definecolor{rouge}{rgb}{1,0,0}
\definecolor{bleu}{rgb}{0,0,1}
\definecolor{vert}{rgb}{0,0.5,0}

\graphicspath{{Figures/}}

\begin{document}

\begin{frontmatter}

\title{Propagation and non-reciprocity in time-modulated diffusion through the lens of high-order homogenization}

%
\author{
M. Touboul$^{1,2}$, B. Lombard$^{3}$, R.C. Assier$^{4}$, S. Guenneau$^{2,5}$ and R.V. Craster$^{1,2,6}$}

\address{$^{1}$Department of Mathematics, Imperial College London, Huxley Building, Queen’s Gate, London SW7 2AZ, UK \\
$^{2}$UMI 2004 Abraham de Moivre-CNRS, Imperial College London, London SW7 2AZ, UK\\
$^{3}$ Aix Marseille Univ, CNRS, Centrale Med, LMA, Marseille, France\\
$^{4}$Department of Mathematics, University of Manchester, Oxford Road, Manchester, M13 9PL, UK \\
$^{5}$The Blackett Laboratory, Department of Physics, Imperial College London, London, SW7 2AZ, UK \\
$^{6}$ Department of Mechanical Engineering, Imperial College London, London SW7 2AZ, UK \\
}

\title{Propagation and non-reciprocity in time-modulated diffusion through the lens of high-order homogenization}




\begin{abstract}
The homogenization procedure developed here is conducted on a laminate with periodic space-time modulation on the fine scale: at leading order, this modulation creates convection in the low-wavelength regime if both parameters are modulated. However, if only one parameter is modulated, which is more realistic, this convective term disappears and one recovers a standard diffusion equation with effective homogeneous parameters; this does not describe the non-reciprocity and the propagation of the field observed from exact dispersion diagrams. This inconsistency is corrected here by considering second-order homogenization which results in a non-reciprocal propagation term that is proved to be non-zero for any laminate and verified via numerical simulation. The same methodology is also applied to the case when the density is modulated in the heat equation, leading therefore to a corrective advective term which cancels out non-reciprocity at the leading order but not at the second order. 

\end{abstract}


\end{frontmatter}



\section{Introduction}

Reciprocity is a fundamental property of classical wave propagation and diffusion that leads to symmetric fields along  opposite directions. 
Breaking reciprocity has received increasing interest in recent years in different physical settings (acoustics, elasticity, electromagnetism, heat transfer, etc.). For waves, it is motivated by the desirable ability to obtain one-way propagation or preferred propagation directions \cite{Rasmussen2021}. For diffusion, applications of non-reciprocity can be found in thermal diodes  \cite{Maldovan2013,Li2004,Li2021} with applications in energy harvesting \cite{Chun2009}, which again motivates approaches to break reciprocity. Different techniques can be used to obtain non-reciprocity, one of which is to periodically modulate physical parameters both in space and time as we shall see in this article.  \\
The spatiotemporal modulation of the parameters in a wavelike fashion to break reciprocity was proposed 60 years ago for polarized electromagnetic waves  \cite{Oliner1961,Simon1960,Cassedy1963,Cassedy1967}; it has since been introduced and studied for acoustic and elastic waves \cite{Nassar2017,Nassar2020}. However, its introduction for the diffusion equation is more recent with  \cite{Torrent2018} giving the first application of these ideas in that context and proposing an effective medium for which the temperature field satisfies a convection-diffusion equation. The effective convective term is responsible for the non-reciprocal thermal properties. However, if only one of the physical parameters of the diffusion equation ( either the capacity or the conductivity) is modulated, this coupling term vanishes and the effective model is reciprocal, which is inconsistent with numerical and experimental observations.\\
Several subsequent works have explored spatiotemporal diffusive metamaterials and have led to an inconsistency. The effective model of \cite{Torrent2018} is also proposed with an additional study of the phase difference of the variations between the thermal conductivity and heat capacity in \cite{Xu2021,OrdonezMiranda2021} and of the transient regime for the second reference in \cite{OrdonezMiranda2021}. In \cite{Xu2022a}, the interest of space-time modulation with respect to the presence of a flow is discussed. Non-reciprocity for modulation of both parameters is illustrated experimentally in \cite{Camacho2020} for the diffusion of electric charges. In \cite{Li2022} it is proved that, for heat transport, the density modulation alters the form of the equation. It is then asserted that reciprocity is recovered, contrary to what was stated in previous works, motivating the study of non-reciprocity while modulating only the thermal conductivity. An answer is given by \cite{Li2022a}, where numerical and experimental evidence of non-reciprocity are provided when only the conductivity is modulated;  these numerical observations are confirmed by \cite{Xu2022}. In the latter, homogenization is also pushed one order further, but does not achieve non-reciprocity when only one parameter is modulated. \\
For acoustics, i.e. wave propagation, second-order homogenization leads to an effective model which always exhibits non-reciprocity even if only one parameter is modulated \cite{Touboul2024} contrary to results from leading-order homogenization \cite{lurie2007,Nassar2017}. Here we adapt the theory of \cite{Touboul2024} to the diffusion equation (see e.g. \cite{shaposhnikova1997asymptotic} for the case of spatial modulation) and we obtain a similar conclusion, thereby fixing the inconsistency between the effective model and the numerical experiments presented in \cite{Xu2022}. This is not a routine extension since the structure of the wave equation and the diffusion equation are intrinsically different, and there is no correspondence even on the dispersion relation because of the time-modulation. To do this, we develop the second-order homogenized model for the time-modulated diffusion equation accounting for propagation and non-reciprocity even if only one parameter is modulated, which was absent in  the effective models proposed so far \cite{Torrent2018,Xu2022}. A similar conclusion holds for the case where a corrective advective term has to be considered due to a modulation of the density: contrary to the conclusion of \cite{Li2022}, we prove that the system is still non-reciprocal if the homogenization procedure is pushed to second order.\\
Moreover, the effective models proposed in \cite{Torrent2018,Xu2022} rely on an expansion of the Floquet-Bloch solution for a material with sinusoidal modulation. An approximation of the dispersion relation at the leading order \cite{Torrent2018} and first-order \cite{Xu2022} is then proposed and, from that, an effective equation for the macroscopic temperature field is identified. This approach does not allow for an arbitrary laminate (while the square-wave modulation is stated in \cite{Xu2022} to be closer to reality and easier to implement  with moving external fields), nor the occurrence of a source term, nor to get the local correctors. This is however possible when considering two-scale asymptotic expansions of the fields \cite{Bensoussan,Bakhvalov1989,cioranescu1999,Conca1997,allaire2022crime,Holmbom2005,Floden2007,Dehamnia2022,Auriault1983,Auriault2009,Pavliotis2008MultiscaleMA}, which is the formalism adopted in the present paper. Our model therefore describes the non-reciprocal effect so far absent through the lens of homogenization, and is valid for any laminate, while accounting for the presence of a source term. \\
In Section \ref{Sec:Pres}, the 1D diffusive laminate with physical properties modulated in time is presented. The moving frame formalism is also detailed, allowing us to apply Floquet-Bloch theory in order to get the exact dispersion relation in the case of a bilayered medium. The two-scale homogenization is performed on the field in Section \ref{Sec:Homog} up to second order. At each order, the limit case of a single modulated parameter is studied. In Section \ref{Sec:Rho_constant}, the case where only the conductivity is modulated is analysed in more detail: the analytical expression of the macroscopic field and the dispersion relation are given, and the latter is compared to the exact dispersion relation in the case of a bilaminate. Section \ref{Sec:Homog_advection} presents the results of a similar approach when the mass density is modulated in the heat equation and an advection term has to be added in the initial governing equation. 

\section{Time-modulated diffusion equation}\label{Sec:Pres}
\begin{figure}[h]
\begin{center}  
    \includegraphics[width=0.7\textwidth]{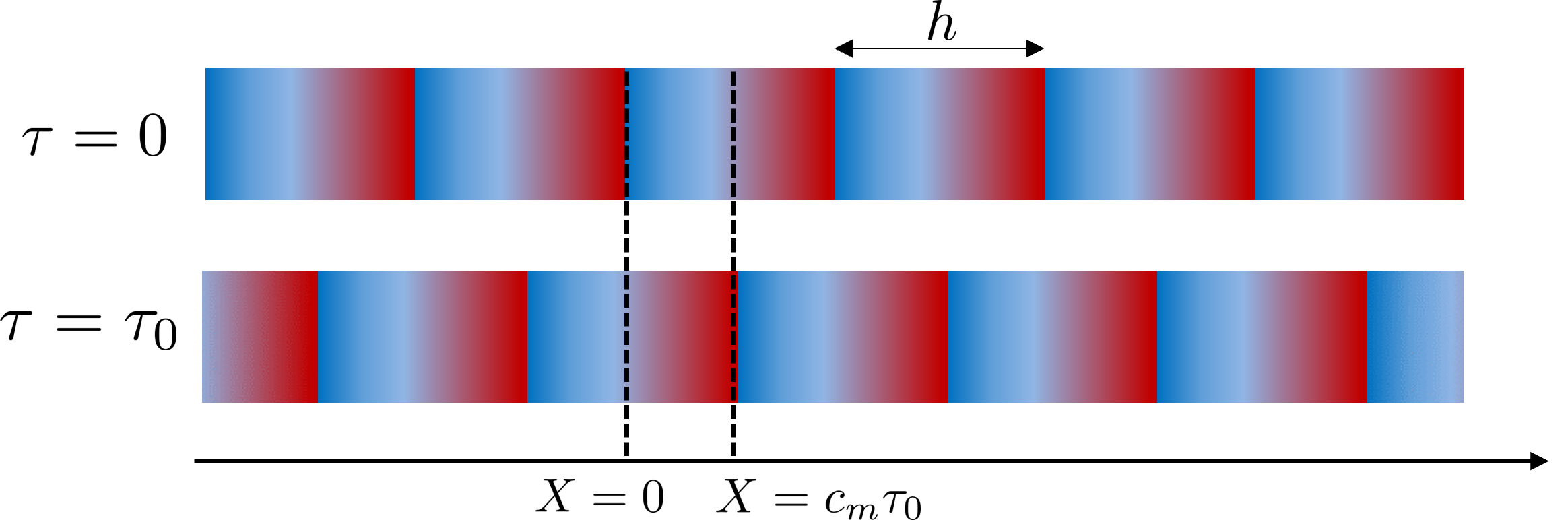}
    \end{center}
    \caption{Bilaminate time-modulated in a wave-like fashion with modulation speed $v_m$ at $\tau=0$ (top) and $\tau=\tau_0$ (bottom).}
        \label{fig:laminate}
\end{figure}

We consider a 1D material whose constitutive parameters $\bar{\sigma}_h>0$ and $\bar{\gamma}_h>0$, corresponding to the conductivity and the capacity of the medium, respectively, are periodically modulated in both space $X$ and time $\tau$ in a wave-like fashion with modulation speed $v_m$:
\begin{equation}
    \label{rho_and_kappa}
    \bar{\sigma}_h(X,\tau) = \bar{\sigma}(X-v_m\tau) \quad \text{ and } \quad \bar{\gamma}_h(X,\tau) = \bar{\gamma}(X-v_m \tau),
\end{equation}
where $\bar{\sigma}$ and $\bar{\gamma}$ are $h$-periodic functions assumed to be piecewise smooth and strictly positive, see Figure \ref{fig:laminate}. \\
Table \ref{table:phys_param} specified the meaning of the field $\Theta_h$ and the constitutive parameters $(\bar{\gamma},\bar{\sigma})$ for two distinct physical settings: thermal diffusion and electric charge diffusion. 
\begin{table}[ht]
\caption{Different meanings of the fields and physical parameters in the diffusion equation, depending on the physical context. }
\begin{center}
\begin{tabular}{ |c|c|c|} 
 \hline
 &Thermal diffusion equation & Electric charge diffusion  \\ \hline 
$\Theta_h$ & Temperature  & Charge density \\  
\hline 
$\bar{{\gamma}}=\bar{\rho}\bar{c}$ & Volumic heat capacity   & Electrical capacity   \\  
& = density $\times$ specific heat capacity & \\  
\hline 
$\bar{\sigma} $ & Thermal conductivity & Electrical conductivity\\  
 \hline
\end{tabular}
\end{center}
\label{table:phys_param}
\end{table}

 Depending on the physical parameters which are modulated in Table \ref{table:phys_param}, different models can be chosen for the time-modulated diffusion equation: the main two are described in the next two subsections. 
\subsection{Model 1: modulation of the capacity and/or conductivity}
A common assumption for the time-modulated diffusion equation is that the energy balance is described by the following equation for the field $\Theta_h$: 
\begin{equation}
    \label{dim_1D_eq}
  \partial_\tau ( \bar{{\gamma}}(\xi) \Theta_h) =  \partial_X\left( \bar{\sigma}(\xi) \partial_X \Theta_h \right)+F
\end{equation}
where $\xi = X-v_m$ and $F(X,\tau)$ is a source term. \\
The temporal modulation is induced by an external field, making the derivation of the time-modulated diffusion equation not straightforward. In the framework of thermoelasticity, it is presented in the Supplementary Material of \cite{Torrent2018}. However, one notes that in \eqref{dim_1D_eq} the capacity is within the time derivative while in \cite{Torrent2018}, the homogenization is performed for a case where the temporal derivative of $\bar{\gamma}$ is cancelled by the external field inducing the modulation. Here, we keep the most general case, without any further assumption on the amplitude of the external field. 
 
\subsection{Model 2: modulation of the density}
In the framework of the heat equation, if the modulation of $\bar{\gamma}$ comes from a modulation of the density $\bar{\rho}$ while the specific capacity $\bar{c}$ is kept constant (see Table \ref{table:phys_param}), a corrective advective term has to be considered so that the conservation of mass holds \cite{Li2022}. In the case of the 1D laminate modulated in a wave-like fashion introduced above, the starting equation \eqref{dim_1D_eq} is then replaced by: 
\begin{equation}
    \label{dim_1D_eq_adv}
 \bar{c}\bar{\rho}(\xi)\partial_\tau \Theta_h+\bar{c}\left[ \bar{\rho}(\xi)-\bar{\rho}_0\right]v_m\partial_X\Theta_h =  \partial_X\left( \bar{\sigma}(\xi) \partial_X \Theta_h \right)+F 
\end{equation}
with $\bar{\rho}_0=\langle \bar{\rho}\rangle_h$ where we have introduced for any function $\bar{g}$
\begin{equation}
    \label{mean}    
\langle \bar{g}\rangle_h=\frac{1}{h}\int_0^h \bar{g}(\xi)\mathrm{d}\xi.
\end{equation}
A similar argument can hold in the framework of electric charge diffusion with the conservation of charges. However, in this case, it is usually implicitly assumed that the modulation comes from a modulation of $\bar{c}$ \cite{Camacho2020}.

\section{\label{sec:level2}Moving frame}
One introduces a moving frame with coordinates $(\xi,\tilde{\tau})=(X-v_m \tau,\tau)$; the field in this frame is denoted by $\tilde{\Theta}_h(\xi,\tilde{\tau})=\Theta_h(X,\tau)$ and satisfies an equation where the coefficients now depend only on the variable $\xi$ and are $h$-periodic in this variable. 
In this moving frame, for the case of Model 1 given by \eqref{dim_1D_eq},  $\tilde{\Theta}_h(\xi,\tilde{\tau})$ satisfies: 
\begin{equation}
    \label{dim_moving_frame}
\partial_\xi\left( \bar{\sigma}(\xi) \partial_\xi\tilde{\Theta}_h(\xi,\tilde{\tau})+v_m\bar{\gamma}(\xi)\tilde{\Theta}_h(\xi,\tilde{\tau})  \right)+\tilde{F}(\xi,\tilde{\tau})=\bar{\gamma}(\xi)\partial_{\tilde{\tau}} \tilde{\Theta}_h(\xi,\tilde{\tau}).
\end{equation}

From \eqref{dim_moving_frame}, we can deduce continuity at any interface between two different layers for both $\tilde{\Theta}_h$ and $ \bar{\sigma} \partial_\xi\tilde{\Theta}_h+v_m\bar{\gamma}\tilde{\Theta}_h$ in the moving frame.
Consequently, we get continuity for ${\Theta}_h$ and $ \bar{\sigma}\partial_X{\Theta}_h+v_m\bar{\gamma}{\Theta}_h $ in the fixed frame across any interface $X_I$.

\subsection{Floquet-Bloch analysis for a bilayered medium}\label{Sec:BF}
Let us use Floquet-Bloch analysis in the moving frame to find the dispersion relation between the Bloch wavenumber $\tilde{\kappa}\in [0,\pi/h]$ which is fixed to be real and the angular frequency $\tilde{\Omega}$. Due to the $h$-periodicity in $\xi$ of the coefficients, we look for solutions in the form
\begin{equation}
    \label{eq:BF_decompo}
    \tilde{\Theta}_h(\xi,\tilde{\tau}) = \mathrm{e}^{\mathrm{i}\tilde{\Omega}\tilde{\tau}}u_h(\xi) \text{ with }u_h(\xi) = \mathrm{e}^{-\mathrm{i}\tilde{\kappa}\xi} \Phi(\xi) 
\end{equation}
with $\Phi$ being $h$-periodic. \\
Using the Floquet-Bloch conditions 
\begin{equation}
    \label{eq:BF_cond}
    \left\lbrace
    \begin{aligned}   
    u_h(h^+) &= \mathrm{e}^{-\mathrm{i}\tilde{\kappa}h}u_h(0^+) \\
   u'_h(h^+) &= \mathrm{e}^{-\mathrm{i}\tilde{\kappa}h}u'_h(0^+), 
   \end{aligned} 
   \right.
\end{equation}
together with the continuity conditions, the dispersion relation can be found as the zeros of  $\mathrm{det}(\boldsymbol{M}(\tilde{\kappa},\tilde{\Omega}))$ for some $4\times4$ matrix $\boldsymbol{M}$. 
In the case of a bilayered material alternating two phases of parameters $(\bar{\sigma}_j,\bar{\gamma}_j)$ ($j=A,B$), and of lengths $\varphi h$ and $(1-\varphi)h$, the continuity conditions in the unit cell $(0,h)$ are written down at $\varphi h$, and at $h$ together with the Floquet-Bloch conditions \eqref{eq:BF_cond}. \\
Because $\kappa X - \Omega \tau = \tilde{\kappa}\xi - \tilde{\Omega}\tilde{\tau}$, the wavenumber $\kappa$ and the angular frequency $\Omega$ in the fixed frame are then given by 
\begin{equation}
    \label{link_frequency_MF_FF}
    \kappa=\tilde{\kappa} \quad \text{ and } \quad \Omega = \tilde{\Omega}+v_m \tilde{\kappa}.
\end{equation}
\subsubsection{Model 1}
For Model 1, the entries of the matrix $\boldsymbol{M}$ are given by: 
\begin{equation}
    \label{eq:Matrix_M}
    \left\lbrace
    \begin{aligned}
   & M_{11} =  \mathrm{e}^{-\mathrm{i}\tilde{\kappa}h}, \qquad  M_{12} =  \mathrm{e}^{-\mathrm{i}\tilde{\kappa}h}, \qquad   M_{13}=- \mathrm{e}^{r_{B}^{+} h}, \qquad M_{14} = - \mathrm{e}^{r_{B}^{-} h}, \\
     & M_{21} =  r_{A}^{+}\bar{\sigma}_A\mathrm{e}^{-\mathrm{i}\tilde{\kappa}h}, \qquad  M_{22} = r_{A}^{-}\bar{\sigma}_A\mathrm{e}^{-\mathrm{i}\tilde{\kappa}h}, \qquad \\
     & M_{23}=- \mathrm{e}^{r_{B}^{+} h}(\bar{\sigma}_Br_{B}^{+}+v_m(\bar{\gamma}_B-\bar{\gamma}_A)),  \qquad M_{24} = - \mathrm{e}^{r_{B}^{-} h}(\bar{\sigma}_Br_{B}^{-}+v_m(\bar{\gamma}_B-\bar{\gamma}_A)), \\
        & M_{31} =  \mathrm{e}^{r_{A}^{+}\varphi h}, \qquad  M_{32} =  \mathrm{e}^{r_{A}^{-}\varphi h}, \qquad   M_{33}=- \mathrm{e}^{r_{B}^{+} \varphi h}, \qquad M_{34} = - \mathrm{e}^{r_{B}^{-} \varphi h}, \\
      & M_{41} = \mathrm{e}^{r_{A}^{+}\varphi h}(\bar{\sigma}_Ar_{A}^{+}+v_m\bar{\gamma}_A) , \qquad  M_{42} =\mathrm{e}^{r_{A}^{-}\varphi h}(\bar{\sigma}_Ar_{A}^{-}+v_m\bar{\gamma}_A), \qquad \\
     & M_{43}=- \mathrm{e}^{r_{B}^{+} \varphi h}(\bar{\sigma}_Br_{B}^{+}+v_m\bar{\gamma}_B),  \qquad M_{44} = - \mathrm{e}^{r_{B}^{-} \varphi h}(\bar{\sigma}_Br_{B}^{-}+v_m\bar{\gamma}_B),\\    
    \end{aligned}
    \right.
\end{equation}
where for $j=A, B$
\begin{equation}
    \label{roots_BF}
    r_{j}^{\pm} = \frac{-v_m \bar{\gamma}_j \pm \sqrt{\Delta_j}}{2\bar{\sigma}_j} \text{ with } {\Delta_j}=v_m^2\bar{\gamma}_j^2+4\bar{\sigma}_j\bar{\gamma}_j\mathrm{i}\Omega.
\end{equation}

\begin{figure}[h]
    \subfloat[Imaginary part \label{fig:DD_both_im}]{%
    \includegraphics[width=0.4\textwidth]{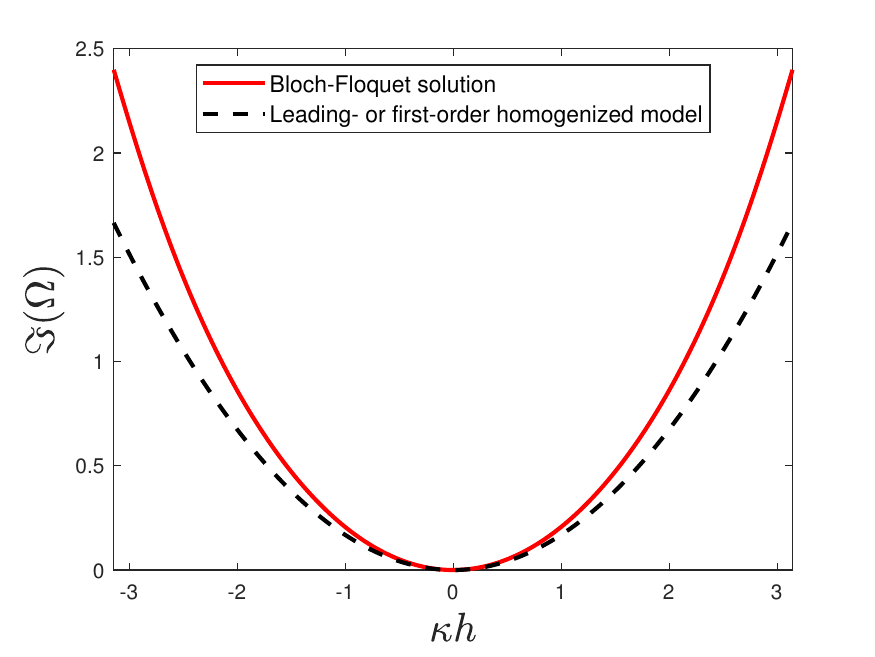}
    }
    \hfill
    \subfloat[Real part \label{fig:DD_both_real}]{%
    \includegraphics[width=0.4\textwidth]{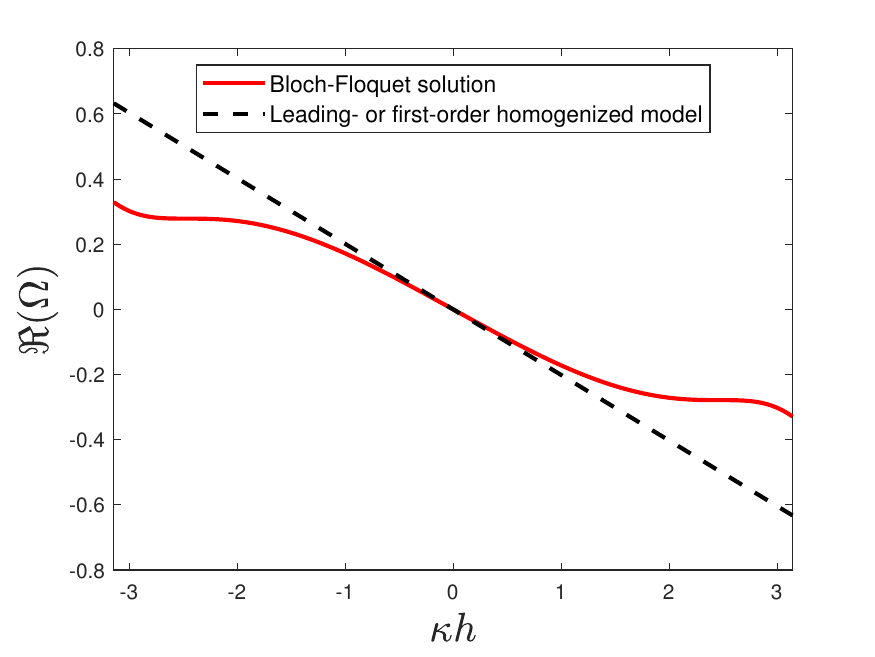}
    }
    \caption{Dispersion diagrams for Model 1 when both parameters are modulated. Exact one in plain lines, leading-order and first-order homogenized model in dashed lines. Non-reciprocity arises from $\Re(\Omega)\neq 0$. The propagation velocity $\Re(\Omega)/\kappa$ being negative, there is propagation against the direction of the space-time modulation.}
    \label{fig:DD_both}
\end{figure}

\begin{figure}[h]
    \subfloat[Imaginary part \label{fig:DD_leading_rho1_im}]{%
    \includegraphics[width=0.4\textwidth]{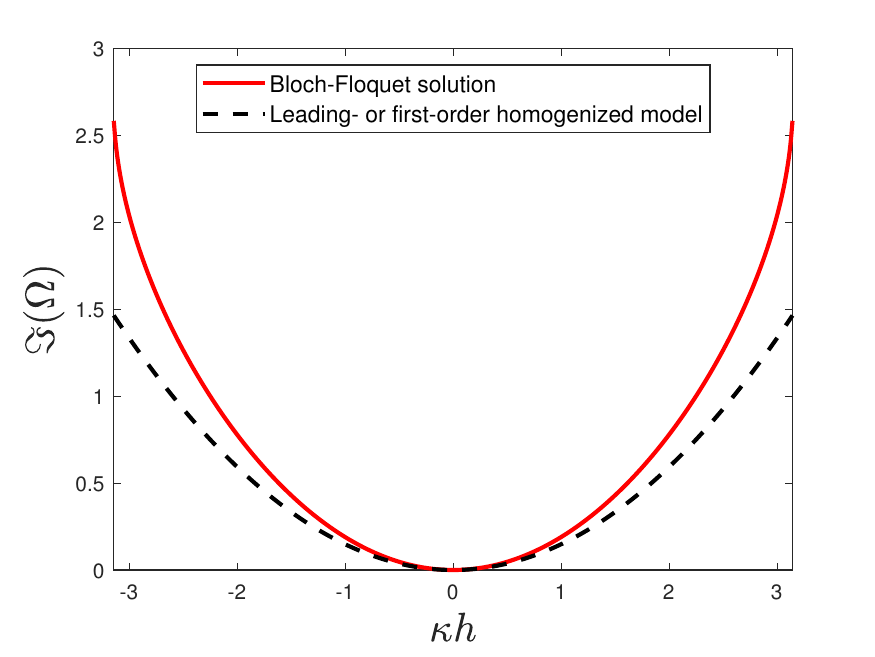}
    }
    \hfill
    \subfloat[Real part \label{fig:DD_leading_rho1_real}]{%
    \includegraphics[width=0.4\textwidth]{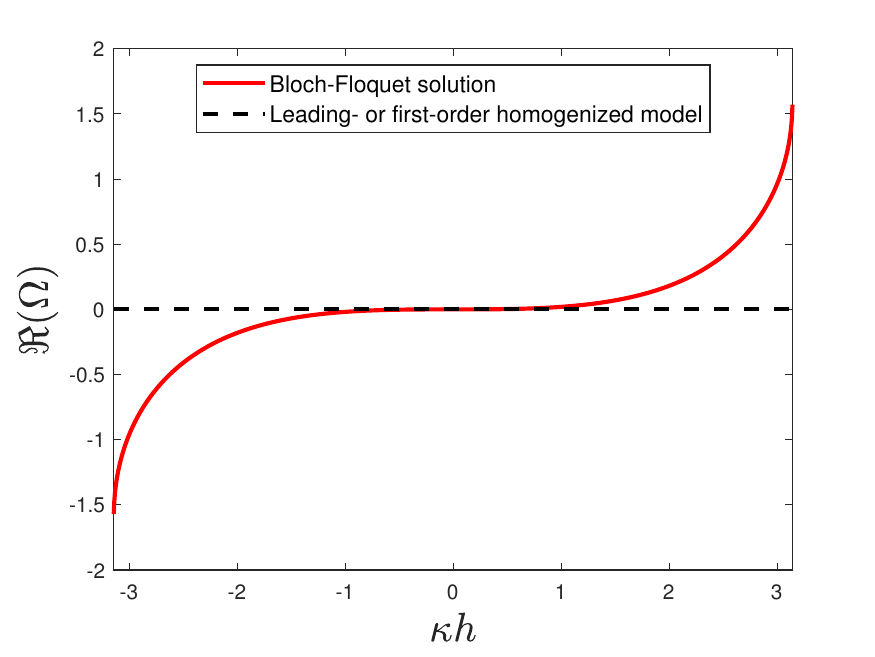}
    }
    \caption{Dispersion diagrams when only $\bar{\sigma}$ is modulated (in this case Model 1 and Model 2 are the same). Exact one in plain lines, leading-order and first-order homogenized model in dashed lines. Non-reciprocity arises from $\Re(\Omega)\neq 0$ but is missed by leading-order models. The propagation velocity $\Re(\Omega)/\kappa$ being positive, there is propagation in the direction of the space-time modulation. }
    \label{fig:DD_leading_rho1}
\end{figure}

The resulting dispersion diagrams are plotted in plain lines in Figure \ref{fig:DD_both} when both parameters are modulated in time with 
\begin{equation}
    \label{eq:bilayer_param}
    \left\lbrace
    \begin{aligned}
        &(\bar{\sigma}_A,\bar{\gamma}_A)=(10\mbox{ kg\,m\,s}^{-3}\,\mbox{K}^{-1},2\times 10^6\mbox{ kg\,m}^{-1}\,\mbox{s}^{-2}\,\mbox{K}^{-1}),\\
        &(\bar{\sigma}_B,\bar{\gamma}_B)=(190\mbox{ kg\,m\,s}^{-3}\,\mbox{K}^{-1},1.4\times 10^6\mbox{ kg\,m}^{-1}\,\mbox{s}^{-2}\,\mbox{K}^{-1}),
    \end{aligned}
    \right.
\end{equation}
a modulation velocity $v_m=5\times 10^{-3}\mbox{ m s}^{-1}$, and a volume fraction $\varphi=0.3$; these parameters are typical of the space-time modulated heat equation considered in \cite{Li2022}.\\
The imaginary part $\Im(\Omega)$ represents the decay rate while the real part $\Re(\Omega)$ reflects the propagative behaviour of the field with a speed $\Re(\Omega)/\kappa$. For non-modulated materials, $\Omega$ is purely imaginary and there is no propagation. On the contrary, as soon as there is a space-time modulation, the field velocity is no longer zero as can be seen in Figure \ref{fig:DD_both_real}. Moreover, the real part of $\Omega$ and $\kappa$ are of opposite sign, therefore leading to non-reciprocity and propagation against the direction of the modulation. \\
The effective models proposed for sinusoidal modulations in \cite{Torrent2018} and \cite{Xu2022} at leading order and first order, respectively,  are extended to the case of any modulated laminate in Section \ref{Sec:Homog}, therefore allowing to deal with the present case of a bilayer (or square modulation). More importantly, as discussed in \cite{Torrent2018} and \cite{Xu2022}, if only one parameter is modulated, the propagative part is completely missed by these effective models. This is confirmed by the dashed lines which stand for the dispersion relation associated to these effective models. We see good agreement around the origin of the diagrams when two parameters are modulated in Figure \ref{fig:DD_both}. However, as we set $\bar{\gamma}_A=\bar{\gamma}_B=2\,10^6\,\mbox{kg\,m}^{-1}\,\mbox{s}^{-2}\,K^{-1}$ in Figure \ref{fig:DD_leading_rho1}, the analytical solution confirms that there is still propagation and non-reciprocity, as soon as one gets further away from the origin. This phenomenon is completely missed by the leading-order and first-order homogenization and will be accounted for in Section \ref{Sec:Homog} through second-order homogenization. 
\subsubsection{Model 2}
\begin{figure}[h]
    \subfloat[Imaginary part \label{fig:DD_leading_rho1_im}]{%
    \includegraphics[width=0.4\textwidth]{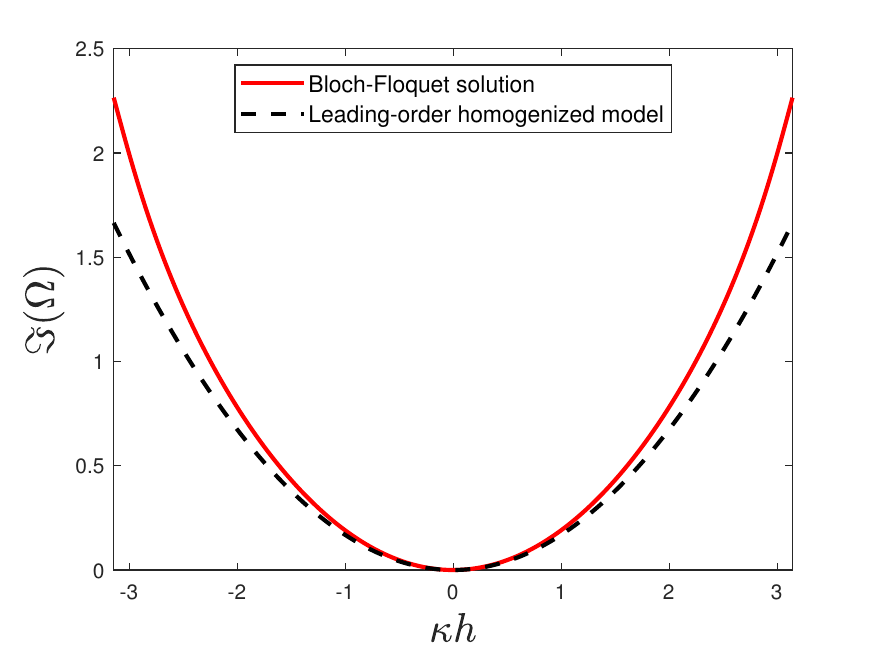}
    }
    \hfill
    \subfloat[Real part \label{fig:DD_leading_rho1_real}]{%
    \includegraphics[width=0.4\textwidth]{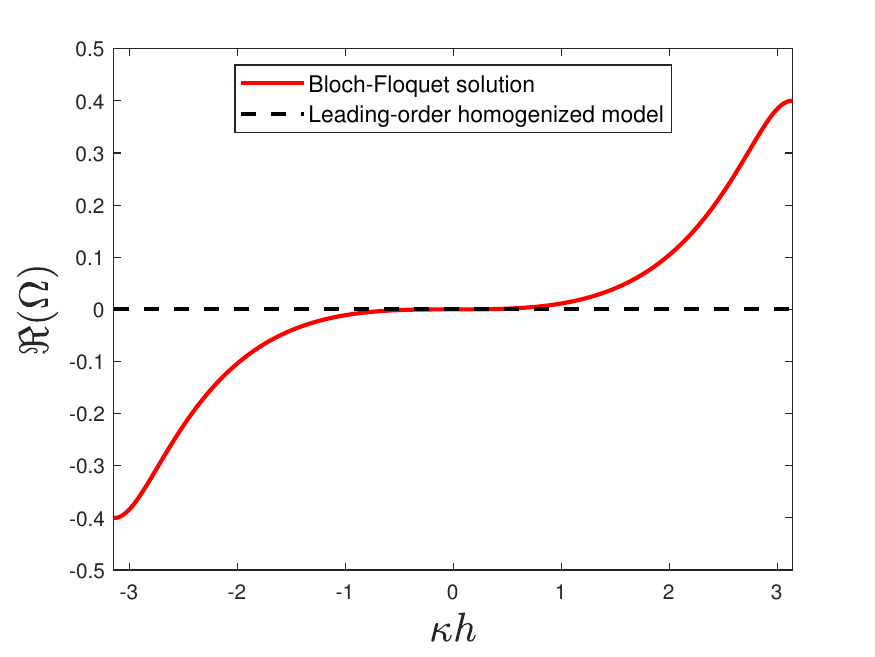}
    }
    \caption{Dispersion diagrams for Model 2 ($\rho$ is modulated requiring therefore a corrective advective term in the modulated diffusion equation). Exact one in plain lines, leading-order homogenized model in dashed lines. Non-reciprocity arises from $\Re(\Omega)$ i.e. the propagation velocity in the direction of the space-time modulation but is missed by leading-order models.}
    \label{fig:DD_leading_adv}
\end{figure}

For Model 2, the entries of the matrix $\boldsymbol{M}$ are given by: 
\begin{equation}
    \label{eq:Matrix_M}
    \left\lbrace
    \begin{aligned}
   & M_{11} =  \mathrm{e}^{-\mathrm{i}\tilde{\kappa}h}, \qquad  M_{12} =  \mathrm{e}^{-\mathrm{i}\tilde{\kappa}h}, \qquad   M_{13}=- \mathrm{e}^{r_{B}^{+} h}, \qquad M_{14} = - \mathrm{e}^{r_{B}^{-} h}, \\
     & M_{21} =  r_{A}^{+}\bar{\sigma}_A\mathrm{e}^{-\mathrm{i}\tilde{\kappa}h}, \qquad  M_{22} = r_{A}^{-}\bar{\sigma}_A\mathrm{e}^{-\mathrm{i}\tilde{\kappa}h}, \qquad  M_{23}=- \bar{\sigma}_Br_{B}^{+}\mathrm{e}^{r_{B}^{+} h},  \qquad M_{24} = - \bar{\sigma}_Br_{B}^{-}\mathrm{e}^{r_{B}^{-} h}, \\
        & M_{31} =  \mathrm{e}^{r_{A}^{+}\varphi h}, \qquad  M_{32} =  \mathrm{e}^{r_{A}^{-}\varphi h}, \qquad   M_{33}=- \mathrm{e}^{r_{B}^{+} \varphi h}, \qquad M_{34} = - \mathrm{e}^{r_{B}^{-} \varphi h}, \\
      & M_{41} = \bar{\sigma}_Ar_{A}^{+}\mathrm{e}^{r_{A}^{+}\varphi h}, \qquad  M_{42} =\bar{\sigma}_Ar_{A}^{-}\mathrm{e}^{r_{A}^{-}\varphi h}, \qquad  M_{43}=-\bar{\sigma}_Br_{B}^{+} \mathrm{e}^{r_{B}^{+} \varphi h},  \qquad M_{44} = - \bar{\sigma}_Br_{B}^{-}\mathrm{e}^{r_{B}^{-} \varphi h},\\    
    \end{aligned}
    \right.
\end{equation}
with for $j=A, B$
\begin{equation}
    \label{roots_BF}
    r_{j}^{\pm} = \frac{-v_m \bar{\rho}_0 \bar{c} \pm \sqrt{\Delta_j}}{2\bar{\sigma}_j} \text{ with } {\Delta_j}=(v_m\bar{\rho}_0 \bar{c} )^2+4\bar{\sigma}_j\bar{\gamma}_j\mathrm{i}\Omega.
\end{equation}
The resulting dispersion diagrams and their leading-order approximations are plotted in Figure \ref{fig:DD_leading_adv} when both parameters are modulated in time following \eqref{eq:bilayer_param}, with $c=1000\,\mbox{m}^2\,\mbox{K}^{-1}\,\mbox{s}^{-2}$, a modulation velocity $v_m=5\cdot 10^{-3}$ $\mathrm{m}\cdot\mathrm{s}^{-1}$, and a volume fraction $\varphi=0.3$ \cite{Li2022}.\\
As for the case of only one parameter modulated in Model 1, we observe that the leading-order homogenized model misses the fact that there is still non-reciprocity and propagation as one gets further away from the origin. This will be corrected by pushing the homogenization up to the second order in Section \ref{Sec:Homog_advection}.

\section{Homogenization for Model 1}\label{Sec:Homog}
\subsection{Non-dimensionalization}
We first introduce some characteristic dimensional wavenumber $\kappa^\star$, conductivity $\sigma^\star$, capacity $\gamma^\star$, field $\Theta^\star$ so that we can define the non-dimensionalized quantities 
\begin{equation}
\begin{aligned}
&\eta = \kappa^\star h, \ x=\kappa^\star X, \ t=(\kappa^\star)^2 \alpha^\star\,\tau,  \sigma=\frac{\bar{\sigma}}{\sigma^\star}, \ \gamma=\frac{\bar{\gamma}}{\gamma^\star}, \ v = \frac{v_m}{v^\star}, \\
&\ k=\frac{\kappa}{\kappa^\star}, \ \omega=\frac{\Omega}{\alpha^\star (\kappa^\star)^2}, T_\eta=\frac{1}{\Theta^\star}\Theta_h, \ f=\frac{1}{\sigma^\star\Theta^\star (\kappa^\star)^2}F,  
  \label{eq:nondimparam}
\end{aligned}
\end{equation}
with $\alpha^\star = \sigma^\star/\gamma^\star$ and $v^\star = \alpha^\star \kappa^\star$. Consequently, one can define the following non-dimensionalized time-modulated diffusion equation
\begin{equation}
    \label{non_dim_1D_eq}
       \partial_t \left({\gamma\left(\frac{x-vt}{\eta}\right)} {T_\eta}\right) =\partial_x\left( {\sigma\left(\frac{x-vt}{\eta}\right)}\partial_x {T_\eta} \right)+{f}
\end{equation}
together with continuity for $T_\eta$ and $\sigma\partial_x T_\eta+v\gamma T_\eta.$
\subsection{Two-scale analysis}
We assume that we have a subwavelength structure, i.e. that the periodicity $h$ is much smaller than the field characteristic wavelength $\lambda^\star = 2\pi/\kappa^\star$. Therefore  the parameter $\eta$ defined in \eqref{eq:nondimparam} is a small parameter: $\eta\ll 1$. Following the two-scale asymptotic technique, the spatial variations of the temperature field $T_\eta$ are assumed to depend on both the macroscopic variable $x$ and the microscopic variable $y=\frac{x-vt}{\eta}$ over the unit periodic cell. The macroscopic variable $x$ describes the slow continuous variations of the field while the microscopic variable $y$ accounts for the small-scale variations, with the parameters $\sigma$ and $\gamma$ varying on this fine scale. Moreover, $x$ and $y$ are assumed to be independent variables, which yields 
\begin{equation}
    \label{eq:diff}
    {\partial_x} \leftrightarrow
{\partial_x} + \frac{1}{\eta}{\partial_y} \text{ and }{\partial_t} \leftrightarrow
{\partial_t} - \frac{v}{\eta} {\partial_y}.
\end{equation}
The field $T_\eta$ is hence expanded using the following ansatz
\begin{equation}
    \label{eq:ansatz_T}
    T_\eta(x,t) = \sum_{j\geq 0} \eta^jT_j(x,y,t).
\end{equation}
Consequently, inserting \eqref{eq:ansatz_T} and \eqref{eq:diff} in \eqref{non_dim_1D_eq}, the non-dimensionalized governing equation reads:
\begin{equation}
    \label{eq:non_dim_two_variables}  
    \left(\partial_x+\frac{1}{\eta}\partial_y \right)\left(\sigma(y) \left(\partial_x+\frac{1}{\eta} \partial_y\right)\sum_{j\geq 0} \eta^jT_j\right)+f(x,t)=    \left(\partial_t-\frac{v}{\eta}\partial_y\right)\left(\gamma(y)\sum_{j\geq 0} \eta^jT_j \right).
\end{equation}
We also define the flux
\begin{equation}
    \label{def_aj_mj}
    a_j = -\sigma(\partial_x T_j+\partial_y T_{j+1}).
\end{equation}
To satisfy the continuity and periodicity conditions, the $T_j$ are assumed to be continuous with respect to the first variable, 1-periodic with respect to the second variable, and both $T_j$ and $-a_j+v\gamma T_j$ are assumed to be continuous with respect to the second variable. \\
We finally introduce the mean value for any function $g(y)$
\begin{equation}
    \label{mean_g}
    \langle g \rangle = \int_0^1 g(y)\mathrm{d}y,
\end{equation}
and we will use extensively that for $g_1$ and $g_2$ continuous and 1-periodic, we have:
$$\left\langle\frac{\mathrm{d}g_1}{\mathrm{d}y}\right\rangle=0 \text{ and } \left\langle g_1\frac{\mathrm{d}g_2}{\mathrm{d}y}\right\rangle = -\left\langle \frac{\mathrm{d}g_1}{\mathrm{d}y}g_2\right\rangle.$$
\subsection{Homogenized equation at leading order}
The homogenized equation at leading order has been obtained for a sinusoidal modulation in \cite{Torrent2018} in the case where heat capacity is outside the time derivative in \eqref{dim_1D_eq}. It is rederived here for the more general case of the present paper. \\
\subsubsection{Derivation in the general case}
Identifying the terms of orders $\eta^{-2}$ in \eqref{eq:non_dim_two_variables}, we get
\begin{equation}
    \label{eq:identif_etam2}
    \partial_y(\sigma\partial_y T_0) = 0
\end{equation}
which leads to 
\begin{equation}
 \partial_y T_0 = \frac{1}{\sigma(y)}q(x,t),
\end{equation}
for some function $q$. Integrating this equation on a unit cell leads to $q(x,t)=0$ and therefore 
\begin{equation}
    T_0(x,y,t)= \mathcal{T}_0(x,t)
\end{equation}
for some function $\mathcal{T}_0(x,t)$.\\
Collecting the terms of order $\eta^{-1}$ in \eqref{eq:non_dim_two_variables}, we then get
\begin{equation}
    \label{eq:identif_etam1}
    \partial_y(\sigma\partial_y T_1)=-v\partial_y(\gamma)T_0-\partial_y(\sigma)\partial_xT_0
\end{equation}
together with continuity for $T_1$ and $-a_0+v\gamma T_0$. We can consequently write $T_1$ as 
\begin{equation}
    \label{form_T1}
    T_1(x,y,t) = \mathcal{T}_1(x,t)+P(y)\partial_x \mathcal{T}_0(x,t)+Q(y) \mathcal{T}_0(x,t),
\end{equation}
for some function $\mathcal{T}_1$, and with $P$ and $Q$ satisfying the cell problems given in \ref{app:cell_leading}. \\
Then, collecting the terms of order $\eta^0$ in \eqref{eq:non_dim_two_variables}, we get 
\begin{equation}
\label{eq:order_eta0}
    \partial_t(\gamma T_0)-v\partial_y(\gamma T_1) = -\partial_x a_0-\partial_y a_1+f,
\end{equation}
with $a_0$ and $a_1$ given by \eqref{def_aj_mj}.
Averaging this equation, using \eqref{form_T1} and introducing the effective parameters 
\begin{equation}
    \label{effective_parameters_zero_order}
  \gamma_0 = \langle \gamma \rangle \text{, } \sigma_0 = \langle \sigma(P'+1)\rangle \text{, and } W_0 = \langle \sigma Q'\rangle, 
\end{equation}
the leading-order effective equation can finally be written 
\begin{equation}
    \label{homog_order0_adim}
    \gamma_0 \partial_t \mathcal{T}_0=\sigma_0 \partial^2_{xx}\mathcal{T}_0+W_0 \partial_x \mathcal{T}_0 + f.
\end{equation}
Consequently, the leading-order homogenized equation is a convection-diffusion equation which involves a convective coefficient $W_0$. One notices that the additional terms in $\partial^2_{xt} \mathcal{T}_0$, which appear in \cite{Torrent2018} and which are ``disregarded because relevant to the transient regime'', do not appear in the setting of this paper at the leading order.
\subsubsection{If only one parameter is modulated}
If only one parameter is modulated in time, i.e. if either $\gamma\equiv 1$ or $\sigma\equiv 1$ then, due to the definition of $W_0$ in \eqref{effective_parameters_zero_order} and \eqref{eq:IPP1}, $W_0=0$ and we recover the usual diffusion equation:
\begin{equation}
    \label{homog_order0_adim_recip}
    \gamma_0 \partial_t \mathcal{T}_0=\sigma_0 \partial^2_{xx}\mathcal{T}_0+f.
\end{equation}
The leading-order field is then purely attenuated in time and reciprocity is recovered. Therefore, when doing leading-order homogenization, it seems that we need modulation of both the conductivity and capacity to have non-reciprocity and propagation, as underlined by \cite{Torrent2018}. We will show that this limitation can be lifted by pushing the homogenization to higher order. 
\subsection{Homogenized equation at first order}
\subsubsection{Derivation in the general case}
Inserting \eqref{form_T1} in \eqref{eq:order_eta0} and using \eqref{homog_order0_adim} and \eqref{eq:simpl_sigma0}, the governing equation for $T_2$ is given by 
\begin{equation}
    \label{syst:T_2}
    \begin{aligned}
    \partial_y(\sigma\partial_y T_2) & = \left[-v\partial_y(\gamma P)-\sigma Q'-\partial_y(\sigma Q)+ W_0+\frac{W_0}{\sigma_0}\partial_y(\sigma P)\right]\partial_x\mathcal{T}_0\\
    &+\left[ \gamma-\gamma_0-\frac{\gamma_0}{\sigma_0}\partial_y(\sigma P)\right]\partial_t\mathcal{T}_0-v\partial_y(\gamma Q)\mathcal{T}_0 -v\partial_y(\gamma)\mathcal{T}_1-\partial_y(\sigma)\partial_x\mathcal{T}_1+\frac{1}{\sigma_0}\partial_y(\sigma P) f
     \end{aligned}
\end{equation}
together with continuity at interfaces of $T_2$ and $-a_1+v\gamma T_1$. We then write $T_2$ as 
\begin{equation}
\begin{aligned}
    \label{form_T2}
    T_2(x,y,t) &= \mathcal{T}_2(x,t)+P(y)\partial_x \mathcal{T}_1(x,t)+Q(y) \mathcal{T}_1(x,t) \\
    &+R(y)\partial_x\mathcal{T}_0(x,t)+S(y)\partial_t\mathcal{T}_0(x,t)+V(y)\mathcal{T}_0(x,t)+A(y)f(x,t)
    \end{aligned}
\end{equation}
with $R$, $S$, $V$ and $A$ satisfying the cell problems given in \ref{app:cell_first}. \\
Collecting the terms of order $\eta$ in \eqref{eq:non_dim_two_variables} leads to 
\begin{equation}
    \label{identif:ordereta}
 \partial_t(\gamma T_1)-v\partial_y(\gamma T_2) = -\partial_x a_1-\partial_y a_2.    
\end{equation}
We then average this equation, insert \eqref{form_T1} and \eqref{form_T2}, and use $\partial_x$\eqref{homog_order0_adim}, \eqref{eq:IPP4} and \eqref{eq:IPP8} to get
\begin{equation}
    \label{effective_eq_order1}
     \gamma_0 \partial_t \mathcal{T}_1-\sigma_0 \partial^2_{xx}\mathcal{T}_1-W_0 \partial_x \mathcal{T}_1 = \mathcal{F}(\mathcal{T}_0)+\mathcal{A}(f)
\end{equation}
where the source terms are defined by 
\begin{equation}
    \label{source_term_F}
    \mathcal{F}(\mathcal{T}_0) = \left[\langle\sigma(Q+R')\rangle-\frac{W_0}{\sigma_0} \langle \sigma P\rangle -\frac{\sigma_0}{\gamma_0}\langle\gamma Q\rangle\right] \partial^2_{xx}\mathcal{T}_0 + \left[\langle\sigma V'\rangle-\frac{W_0}{\gamma_0}\langle\gamma Q\rangle \right] \partial_x\mathcal{T}_0
\end{equation}
and 
\begin{equation}
    \label{source_term_A}
    \mathcal{A}(f)  = -\frac{1}{\gamma_0}\langle \gamma Q \rangle f.
\end{equation}
\subsubsection{If only one parameter is modulated}
If only the capacity is modulated, i.e. $\sigma\equiv 1$,  one has $P\equiv 0$, $W_0\equiv 0$, and the source terms write $$\mathcal{F}(\mathcal{T}_0)+\mathcal{A}(f)=-\frac{\langle\gamma Q\rangle}{\gamma_0}(\partial^2_{xx}\mathcal{T}_0+f).$$
Moreover, integrating \eqref{cell_pb_Q}, we can show that 
\begin{equation}
    \label{eq:rhoQ_zero}
  \langle\gamma Q\rangle=\langle\gamma\rangle\langle Q\rangle-\frac{1}{c}\langle QQ'\rangle=0,  
\end{equation}
the last identity coming from $\langle Q\rangle=0$ and from $\langle QQ'\rangle =\frac{1}{2}\langle (Q^2)'\rangle$ together with periodicity conditions.
Consequently, the first-order effective equation \eqref{effective_eq_order1}  reduces to 
\begin{equation}
    \label{effective_eq_order1_sigma1}
   \gamma_0\partial_t \mathcal{T}_1-\partial^2_{xx}\mathcal{T}_1 = 0.
\end{equation}
Conversely, if only the conductivity is modulated, i.e. $\gamma\equiv1$, then $Q\equiv 0$, $W_0\equiv 0$, and using \eqref{eq:IPP2} and \eqref{eq:IPP6}, the first-order effective equation \eqref{effective_eq_order1} reads 
\begin{equation}
    \label{effective_eq_order1_rho1}
   \partial_t \mathcal{T}_1-\sigma_0 \partial^2_{xx}\mathcal{T}_1 =0.\end{equation}
Consequently, at first order the same analysis holds than at leading order: if only one parameter is modulated, the convective term vanishes and the non-reciprocity is missed by the low-frequency homogenization as noticed in \cite{Xu2022} for sinusoidal modulation. 
\subsection{Homogenized equation at second order}
\subsubsection{Derivation in the general case}
Inserting \eqref{form_T1} and \eqref{form_T2} in \eqref{identif:ordereta}, using \eqref{homog_order0_adim}, $\partial_x$\eqref{homog_order0_adim} and \eqref{effective_eq_order1}, we get the following governing equation for $T_3$
\begin{equation}
    \label{syst:T_3}
    \begin{aligned}
    \partial_y(\sigma\partial_y T_3) & = \left[-\sigma(Q+R')-\partial_y(\sigma R)+\frac{W_0}{\sigma_0}\sigma P+\frac{\sigma_0}{\gamma_0}\gamma Q-\frac{\sigma_0}{\gamma_0}v\partial_y(\gamma S)-\mathcal{C}_1(\sigma_0+\partial_y(\sigma P))
    \right] \partial^2_{xx}\mathcal{T}_0 \\
    &+\left[-\sigma S'-\partial_y(\sigma S)+\gamma P-\frac{\gamma_0}{\sigma_0}\sigma P\right]\partial^2_{tx}\mathcal{T}_0 \\
    &+\left[ -v\partial_y(\gamma R)-\sigma V'-\partial_y(\sigma V)+\frac{W_0}{\gamma_0}\gamma Q-\frac{W_0}{\gamma_0}v\partial_y(\gamma S)-\mathcal{C}_2(\sigma_0+\partial_y(\sigma P)) \right]\partial_x\mathcal{T}_0 \\
    &-v\partial_y(\gamma V)\mathcal{T}_0+\left[-v\partial_y(\gamma P)-\sigma Q'-\partial_y(\sigma Q)+W_0+\frac{W_0}{\sigma_0}\partial_y(\sigma P)\right]\partial_x\mathcal{T}_1\\
    & +\left[ \gamma-\gamma_0-\frac{\gamma_0}{\sigma_0}\partial_y(\sigma P)\right]\partial_t\mathcal{T}_1-v\partial_y(\gamma Q)\mathcal{T}_1 -v\partial_y(\gamma)\mathcal{T}_2-\partial_y(\sigma)\partial_x\mathcal{T}_2\\
    &+ \left[-v\partial_y(\gamma A)+\frac{\gamma}{\gamma_0}Q-\frac{v}{\gamma_0}\partial_y(\gamma S)-\frac{1}{\gamma_0}\langle\gamma Q\rangle\left(1+\frac{1}{\sigma_0}\partial_y(\sigma P)\right) \right]f\\
    &+\left[-\sigma A'-\partial_y(\sigma A)+\frac{1}{\sigma_0}\sigma P\right] \partial_x f,
     \end{aligned}
\end{equation}
together with continuity of $T_3$ and $-a_2+v\gamma T_2$ at interfaces. The coefficients $\mathcal{C}_1$ and $\mathcal{C}_2$ are defined by
\begin{equation}
    \label{def_C_i}
    \left\lbrace
    \begin{aligned}
    &\mathcal{C}_1= -\frac{1}{\sigma_0}\langle\sigma(Q+R')\rangle+\frac{W_0}{\sigma_0^2}\langle\sigma P\rangle +\frac{1}{\gamma_0}\langle\gamma Q\rangle \\
    & \mathcal{C}_2 = -\frac{v}{\sigma_0}\langle \gamma Q P'\rangle +\frac{W_0}{\gamma_0\sigma_0}\langle\gamma Q\rangle,
    \end{aligned}
    \right.
\end{equation}
where we used \eqref{eq:IPP6} to get the expression of $\mathcal{C}_2$.
Consequently, $T_3$ can be written as 
\begin{equation}
\begin{aligned}
    \label{form_T3}
    T_3(x,y,t) &= \mathcal{T}_3(x,t)+P(y)\partial_x \mathcal{T}_2(x,t)+Q(y) \mathcal{T}_2(x,t)\\
    &+R(y)\partial_x\mathcal{T}_1(x,t)+S(y)\partial_t\mathcal{T}_1(x,t)+V(y)\mathcal{T}_1(x,t)\\
    &+L(y)\partial^2_{xx}\mathcal{T}_0(x,t)+M(y)\partial^2_{tx}\mathcal{T}_0(x,t)+N(y)\partial_x\mathcal{T}_0(x,t) + O(y) \mathcal{T}_0(x,t)\\
    &+B(y)f(x,t)+C(y)\partial_x f(x,t),
    \end{aligned}
\end{equation}
with $L$, $M$, $N$, $O$, $B$ and $C$ solutions of the cell problems given in \ref{app:cell_second}.\\
Collecting the terms of order $\eta^2$ in \eqref{eq:non_dim_two_variables} and averaging on a unit cell, we get 
\begin{equation}
    \label{eq:ident_eta2}
    \partial_t \langle \gamma T_2\rangle = \partial_x\langle\sigma(\partial_x T_2+\partial_y T_3)\rangle.
\end{equation}
Inserting \eqref{form_T2} and \eqref{form_T3} in this equation, and using \eqref{eq:IPP4}, \eqref{effective_eq_order1}, and $\partial_x$\eqref{effective_eq_order1}, we get the final effective equation at the second order
\begin{equation}
    \label{effective_eq_order2}
     \gamma_0 \partial_t \mathcal{T}_2-\sigma_0 \partial^2_{xx}\mathcal{T}_2-W_0 \partial_x \mathcal{T}_2 = \mathcal{F}(\mathcal{T}_1)+\mathcal{E}(\mathcal{T}_0)+\mathcal{B}(f)
\end{equation}
where the source terms are defined by 
\begin{equation}
    \label{source_term_E}
    \begin{aligned}
    \mathcal{E}(\mathcal{T}_0) &= \left[\langle\sigma R\rangle+\langle\sigma L'\rangle+\langle \sigma P\rangle\mathcal{C}_1\right]\partial^3_{xxx}\mathcal{T}_0  + \left[\langle\sigma S\rangle+\langle\sigma M'\rangle \right] \partial^3_{txx}\mathcal{T}_0-\langle\gamma R\rangle\partial_{tx}^2\mathcal{T}_0-\langle\gamma S\rangle\partial_{tt}^2\mathcal{T}_0\\
    & + \left[\langle\sigma V\rangle+\langle\sigma N'\rangle +\langle\gamma Q\rangle\frac{\sigma_0}{\gamma_0}\mathcal{C}_1+\langle\sigma P\rangle\mathcal{C}_2\right] \partial^2_{xx}\mathcal{T}_0 +    \left[\langle\sigma O'\rangle +\langle\gamma Q\rangle\frac{\sigma_0}{\gamma_0}\mathcal{C}_2\right] \partial_{x}\mathcal{T}_0-\langle\gamma V\rangle\partial_t\mathcal{T}_0
    \end{aligned}
\end{equation}
and 
\begin{equation}
\begin{aligned}
    \label{source_term_B}
    \mathcal{B}(f)  &= \left[\langle\sigma A\rangle+\langle \sigma C'\rangle \right] \partial_{xx}^2 f + \left[ \langle\sigma B'\rangle+\frac{\langle\sigma P\rangle\langle\gamma Q\rangle}{\gamma_0\sigma_0}\right]\partial_xf -\langle\gamma A\rangle\partial_tf+\frac{1}{\gamma_0^2}\langle\gamma Q\rangle^2f,
    \end{aligned}
\end{equation}
and $\mathcal{F}$ has been defined in \eqref{source_term_F}.
\subsubsection{If only one parameter is modulated}
If $\sigma \equiv 1$, then $A \equiv 0$, and $\mathcal{C}_1=\mathcal{C}_2=0$ with \eqref{eq:rhoQ_zero}. Therefore, the second-order effective equation \eqref{effective_eq_order2} reads
\begin{equation}
    \label{effective_eq_order2_sigma1}
    \gamma_0 \partial_t \mathcal{T}_2- \partial^2_{xx}\mathcal{T}_2 = \mathcal{E}(\mathcal{T}_0)
\end{equation}
where the source term $\mathcal{E}$ reduces to
\begin{equation}
    \label{eq:source_E_sigma1}
    \mathcal{E}(\mathcal{T}_0) = -\langle\gamma R\rangle \partial^2_{xt}\mathcal{T}_0-\langle\gamma S\rangle\partial^2_{tt}\mathcal{T}_0-\langle\gamma V\rangle\partial_t\mathcal{T}_0.
\end{equation}
Non-reciprocity happens when the parity $x\leftrightarrow -x$ is broken by odd-order spatial derivatives. Consequently, non-reciprocity holds as soon as 
\begin{equation}
    \label{def_NR_term_sigma1}
    \mathcal{N}_\gamma=-\langle\gamma R\rangle
\end{equation}
 is non-zero. Integrating \eqref{cell_pb_R} and using $\langle R'\rangle=\langle Q\rangle = 0$ we get 
\begin{equation}
    \label{eq:Rprim_sigma1}
    R'=-2Q.
\end{equation}
Moreover, integrating \eqref{cell_pb_Q} together with $\langle Q'\rangle = 0$, one gets
\begin{equation}
    \label{eq:rhoQprim_sigma1}
    \gamma = -\frac{1}{v}Q'+\gamma_0.
\end{equation}
Consequently, using \eqref{eq:rhoQprim_sigma1},  $\langle R\rangle=0$ and \eqref{eq:Rprim_sigma1}, $\mathcal{N}_\gamma$ reads
\begin{equation}
    \mathcal{N}_\gamma = \frac{1}{v} \langle Q' R\rangle = -\frac{1}{v} \langle Q R'\rangle = \frac{2}{v} \langle Q^2 \rangle.
\end{equation}
From \eqref{eq:rhoQprim_sigma1}, one can notice that $Q'=v(\gamma_0-\gamma)$, and is therefore not identically zero as soon as $v\neq 0$ and $\gamma$ is not constant.
One can then deduce that $Q\not\equiv 0 $ and 
\begin{equation}
    \mathcal{N}_\gamma  \neq 0.
\end{equation}
In a similar fashion, if $\gamma\equiv 1$, $V=O\equiv 0$, and the effective equation reduces to 
\begin{equation}
    \label{effective_eq_order2_rho1}
    \partial_t \mathcal{T}_2-\sigma_0 \partial^2_{xx}\mathcal{T}_2 = \mathcal{B}(f)+\mathcal{E}(\mathcal{T}_0)
\end{equation}
where the second source term reads: 
\begin{equation}
    \label{eq:source_term_E_rho1}
    \mathcal{E}(\mathcal{T}_0)=\mathcal{N}_\sigma\partial^3_{xxx}\mathcal{T}_0+\left[\langle\sigma S\rangle+\langle\sigma M'\rangle \right]\partial^3_{txx}\mathcal{T}_0+\langle\sigma N'\rangle\partial^2_{xx}\mathcal{T}_0,
\end{equation}
with 
\begin{equation}
    \label{def_NR_term_rho1}
    \mathcal{N}_\sigma=\langle\sigma R\rangle+\langle\sigma L'\rangle-\frac{1}{\sigma_0}\langle\sigma P\rangle\langle\sigma R'\rangle
\end{equation}
the term associated to non-reciprocity. 
In this case, integration by parts of $\langle\eqref{cell_pb_L}\times P - \eqref{cell_pb_P}\times L\rangle$ leads to
\begin{equation}
    \label{eq:IPPLrho1}
    \langle \sigma L' \rangle =\langle\sigma RP'\rangle-\langle\sigma R'P\rangle-\sigma_0 v \langle S'P\rangle+\frac{1}{\sigma_0}\langle\sigma R'\rangle\langle\sigma P\rangle.
\end{equation} 
Moreover, integration of \eqref{cell_pb_R}, and multiplication by $P$ together with $\langle P\rangle =0$ leads to 
\begin{equation}
    \label{sigmaRprim_rho1}
    \langle\sigma R'P\rangle = -v\langle P^2\rangle.
\end{equation}
In a similar way, integration of \eqref{cell_pb_P} and multiplication by $R$ together with $\langle R\rangle =0$ leads to 
\begin{equation}
    \label{sigmaPprim_rho1}
    \langle\sigma P'R\rangle = -\langle \sigma R\rangle. 
\end{equation}
Finally, integration of \eqref{cell_pb_S} together with $\langle S'\rangle =0$ yields 
\begin{equation}
    \label{eq:Sprim_rho1}
    S'=-\frac{1}{\sigma_0}P
\end{equation}
and consequently
\begin{equation}
    \label{sigmaSprim_rho1}
    \langle S'P\rangle = -\frac{1}{\sigma_0}\langle P^2\rangle.
\end{equation}
Inserting \eqref{sigmaPprim_rho1}, \eqref{sigmaRprim_rho1} and \eqref{sigmaSprim_rho1} in \eqref{eq:IPPLrho1} allows to write
\begin{equation}
    \label{eq:nonnulNbeta}
    \mathcal{N}_\sigma = 2v\langle P^2\rangle > 0
\end{equation}
because $P\not\equiv 0$ (otherwise, one would get $\sigma$ constant from \eqref{cell_pb_P}). \\
Consequently, we proved that if only one parameter is modulated, either $\gamma$ or $\sigma$, we still get non-reciprocity through the lens of low-frequency homogenization. We also encapsulate through these odd-order space derivatives terms the propagative part, as will be illustrated later on with dispersion diagrams and time-domain simulations.   
  
\section{Numerical illustration of the effective behaviour for Model 1}\label{Sec:Rho_constant}
\subsection{Effective behaviour at the leading order}
\begin{figure}[h]
    \subfloat[Both parameters modulated (Model 1)\label{fig:field_leading_both}]{%
    \includegraphics[width=0.48\textwidth]{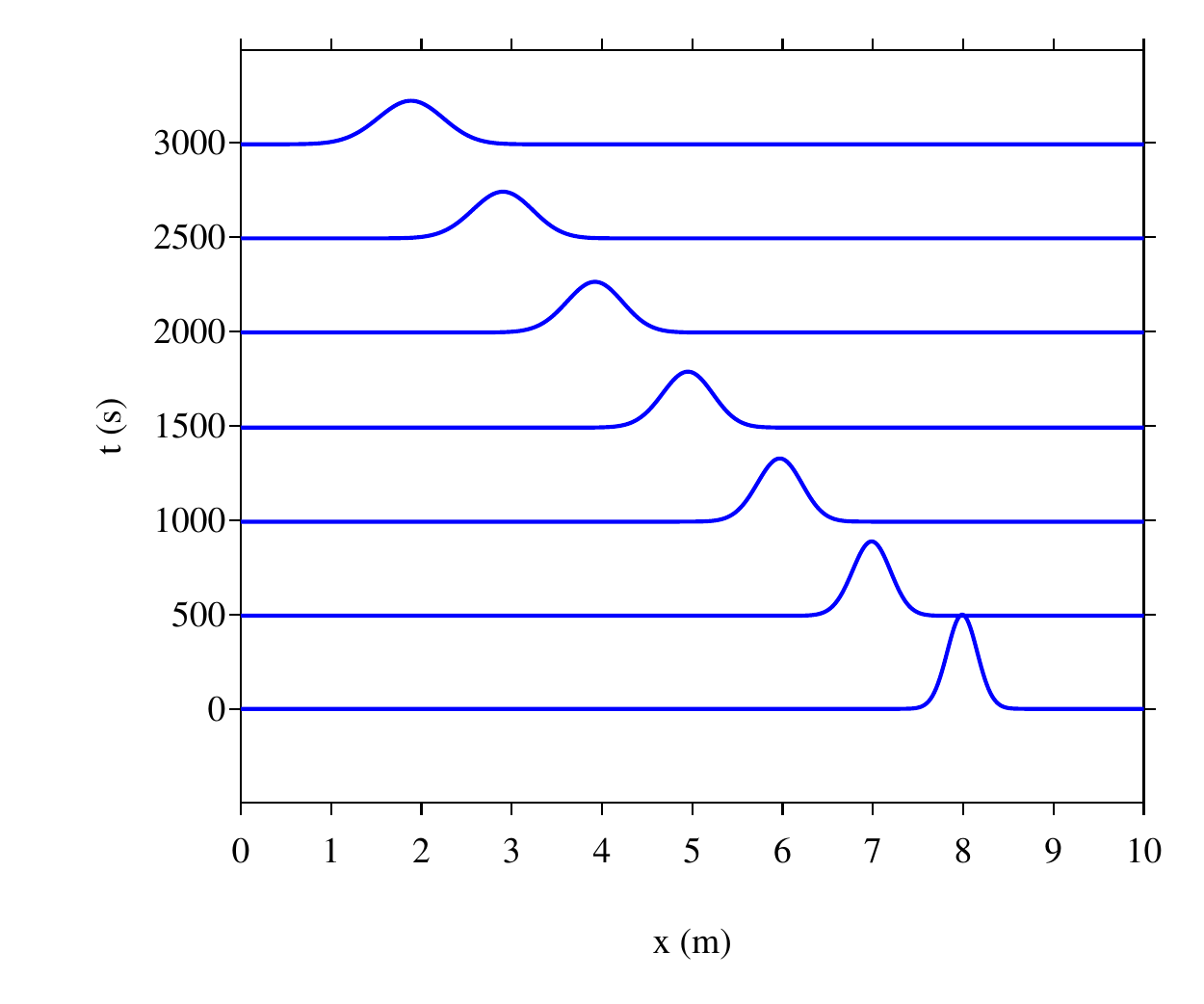}
    }
    \hfill
    \subfloat[$\gamma$ constant \label{fig:field_leading_one}]{%
    \includegraphics[width=0.48\textwidth]{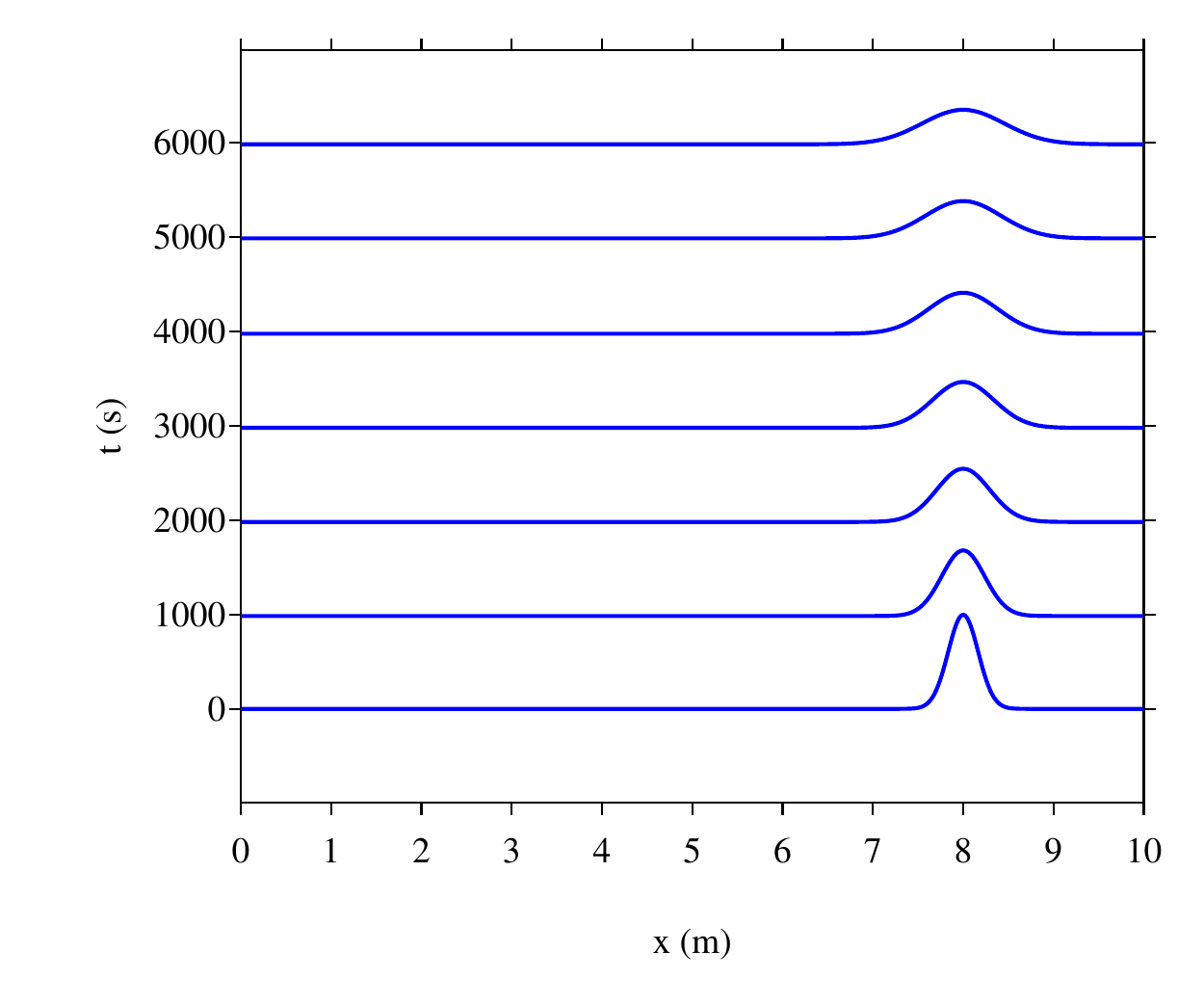}
    }
    \caption{Snapshots of the leading-order field $\Theta_0$ at different times.}
    \label{fig:field_leading_order}
\end{figure}

At leading order, the effective convection-diffusion equation \eqref{homog_order0_adim} reads in dimensional coordinates: 
\begin{equation}
    \label{homog_order0_dim}
    \bar{\gamma}_0 \partial_\tau \mathcal{\Theta}_0-v_m \mathcal{W}_0 \partial_X \mathcal{\Theta}_0 = \bar{\sigma}_0  \partial^2_{XX}\mathcal{\Theta}_0+F
\end{equation}
where the dimensionalized effective parameters are defined by 
\begin{equation}
    \label{dim_eff_param}
    \left\lbrace
    \begin{aligned}
    &\bar{\gamma}_0 = \langle \gamma \rangle_h \\
    &\bar{\sigma}_0 = \left\langle \frac{1}{\bar{\sigma}} \right\rangle_h^{-1} \\
    & \mathcal{W}_0 = \left\langle \frac{\bar{\gamma}}{\bar{\sigma}} \right\rangle_h\bar{\sigma}_0-\bar{\gamma}_0.
    \end{aligned}
    \right.
\end{equation}
For the parameters of Figure \ref{fig:DD_both}, where both parameters are modulated, these effective parameters are given by 
\begin{equation}
    \label{eq:eff_param_both}
    \left\lbrace
    \begin{aligned}
       & \bar{\gamma_0}  = 1.76\times 10^6\,\mbox{kg m}^{-1}\,\mbox{s}\\^{-2}\,\mbox{K}^{-1}
       &\bar{\sigma_0} = 29.7\mbox{kg m}\,\mbox{s}\\^{-3}\,\mbox{K}^{-1} \\
       & v_m \mathcal{W}_0 = 3.54\times 10^3\,\mbox{ kg s}^{-3}\,\mbox{ K}^{-1}.
    \end{aligned}
    \right.
\end{equation}
The Reynolds number $\frac{v_m \mathcal{W}_0}{\bar{\sigma}_0}$ for this effective equation is about 119 $\gg 1$, so front propagation is clearly expected. \\
For numerical examples, we will consider initial data given by 
a Gaussian centered at $X_0$ with standard deviation $\nu$: 
\begin{equation}
\Theta_0(X,0)=\theta_0(X)=\mathrm{e}^{-\left(\frac{X-X_0}{\nu}\right)^2}.
\label{GaussianInit}
\end{equation}
The time and spatial evolution of $\Theta_0$ for the effective equation defined by \eqref{homog_order0_dim} and \eqref{eq:eff_param_both} is computed by finite differences using a  semi-implicit Crank-Nicolson scheme \cite{Strikwerda2004}. The snapshots at different time steps are given in Figure \ref{fig:field_leading_both}, for $X_0 = 8$m, $\nu=1$m, which confirms a clear front propagation and consequently a strong non-reciprocal effective behaviour. \\
However, if only one parameter is modulated, $W_0$ vanishes, and one recovers a symmetric solution, as observed in Figure \ref{fig:field_leading_one} using the same parameters than in the dispersion diagrams presented in Figure \ref{fig:DD_leading_rho1}. 

\subsection{Effective macroscopic field for $\gamma$ constant}
\begin{figure}[h!]
    \subfloat[$\nu = 1  \,(\eta=0.105)$ \label{fig:Homog2_rho1_field_a}]{%
    \includegraphics[width=0.48\textwidth]{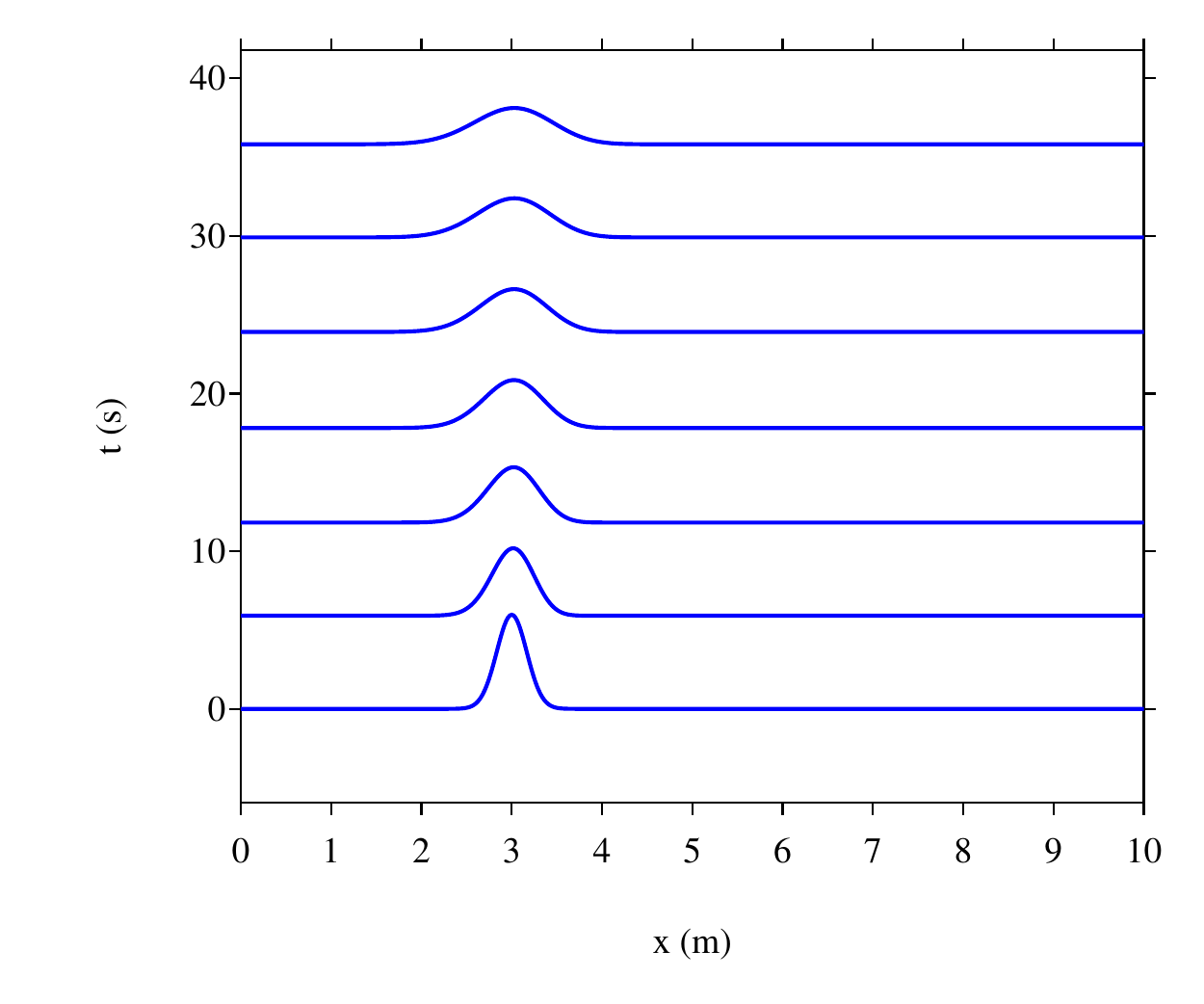}
    }
    \hfill
    \subfloat[$\nu = 0.1  \,(\eta=1.05)$ \label{fig:Homog2_rho1_field_b}]{%
    \includegraphics[width=0.48\textwidth]{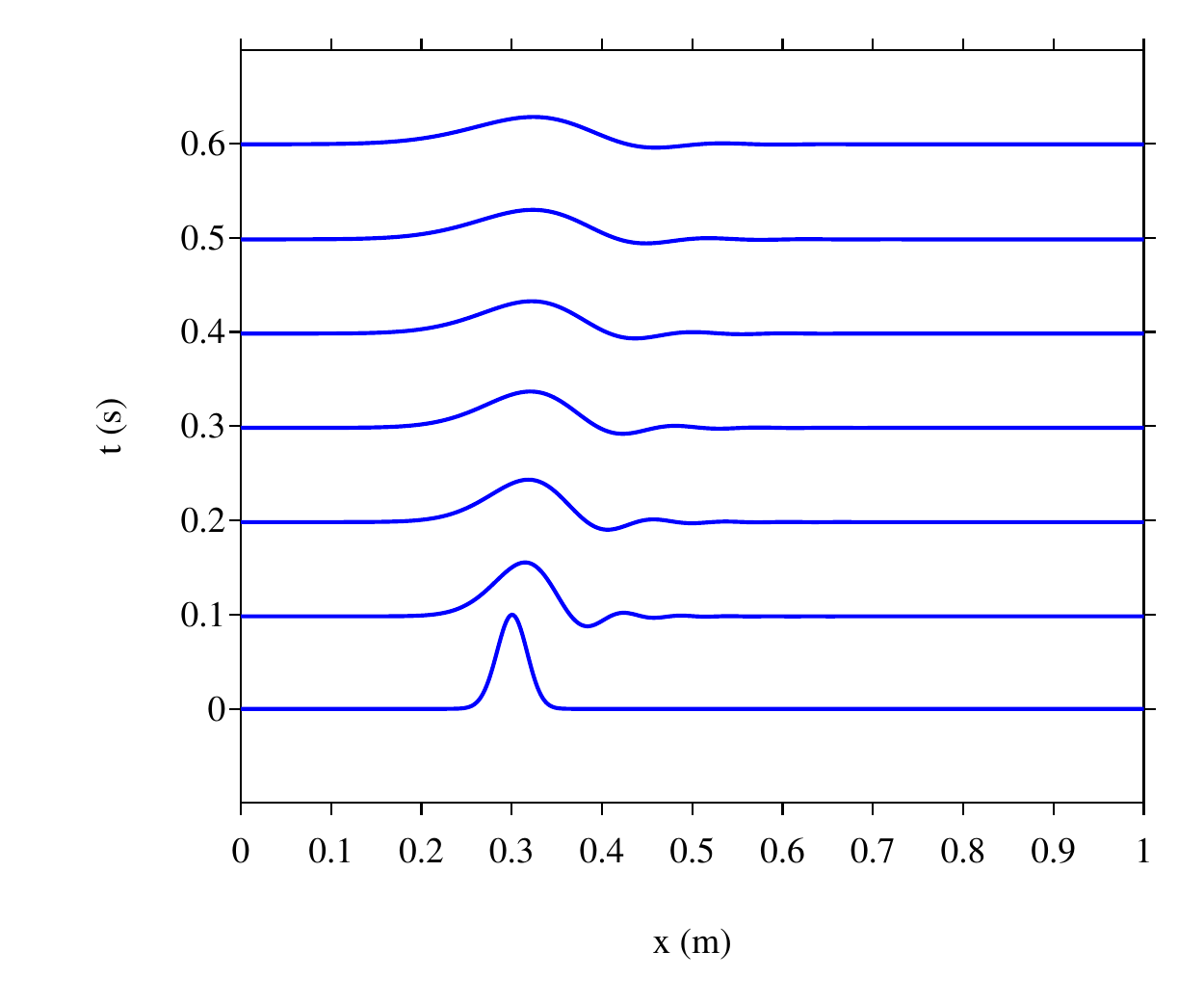}
    }\\
    \centering
    \subfloat[$\nu = 0.05 \,(\eta=2.09)$\label{fig:Homog2_rho1_field_c}]{%
    \includegraphics[width=0.48\textwidth]{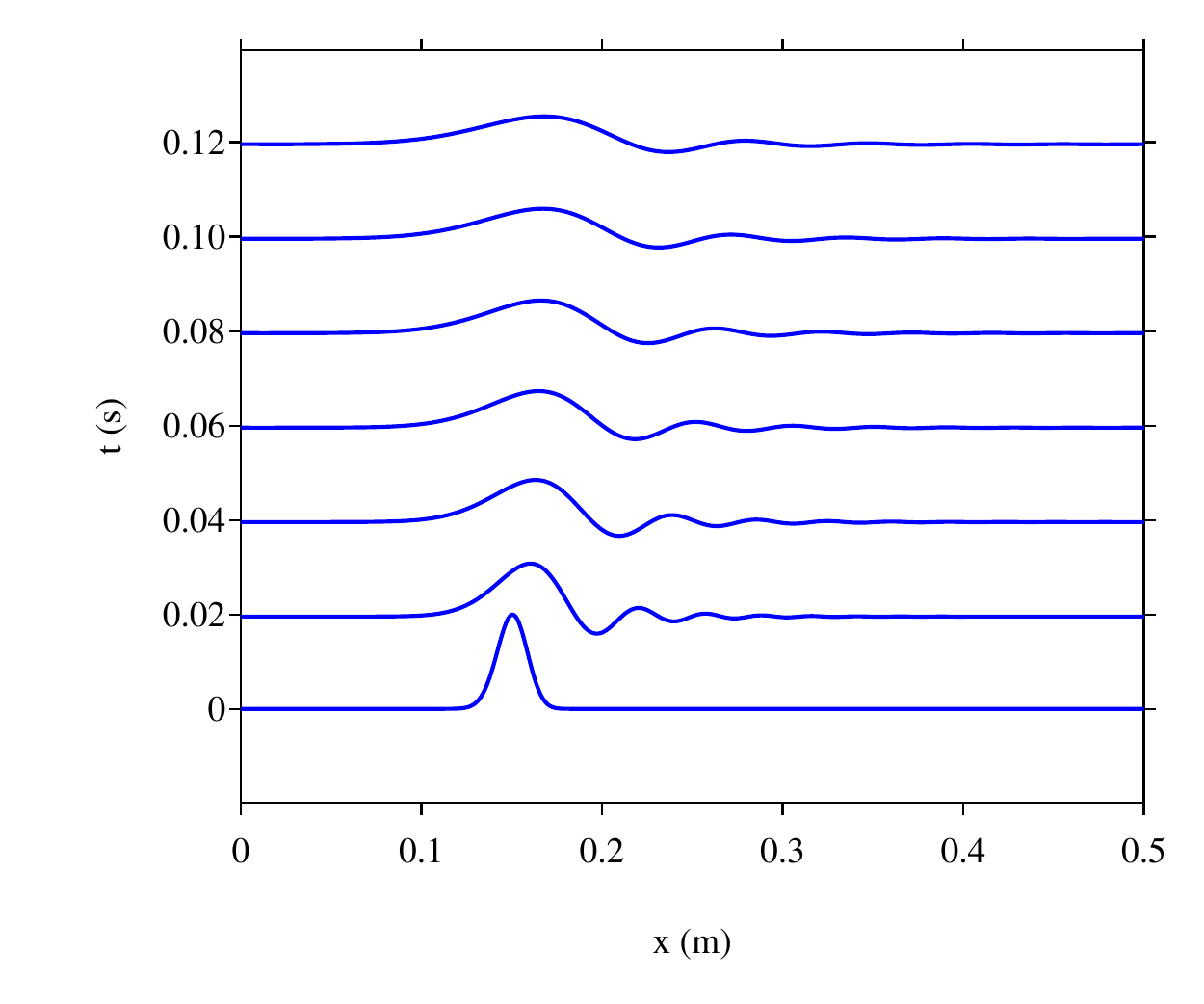}
    }

    \caption{$\gamma$ constant. Snapshots of the second-order homogenized mean field $\Theta^{(2)}$ at different times and for different values of the small parameter $\eta$.}
    \label{fig:Homog2_rho1_field}
\end{figure}

Since it seems more realistic to consider the time-modulated heat equation \eqref{dim_1D_eq} if only the thermal conductivity is modulated \cite{Li2022}, we focus here on the case where the heat capacity $\gamma$ is constant to explicit the macroscopic effective behaviour associated to the second-order effective model developed in Section \ref{Sec:Homog}. However, similar expressions and illustrations can be obtained in the case $\sigma$ constant. \\
First we introduce the second-order macroscopic mean field 
\begin{equation}
    \label{def_mean_field_order2}
    \mathcal{T}^{(2)}=\mathcal{T}_0+\eta\mathcal{T}_1+\eta^2\mathcal{T}_2
\end{equation}
that describes the macroscopic behaviour of the second-order field $T_2$. Combining \eqref{homog_order0_adim_recip}, \eqref{effective_eq_order1_rho1}, \eqref{effective_eq_order2_rho1} this mean field satisfies when $f=0$
\begin{equation}
\label{eq:def_eff_0}
    \partial_t \mathcal{T}^{(2)}-\sigma_0 \partial^2_{xx}\mathcal{T}^{(2)} = \eta^2\mathcal{E}(\mathcal{T}_0)
\end{equation}
with $\mathcal{E}(\mathcal{T}_0)$ given by \eqref{eq:source_term_E_rho1} and \eqref{eq:nonnulNbeta}. We will now write down the coefficients of this equation explicitly in terms of $\sigma$ only. \\
Firstly, integrating \eqref{cell_pb_P} together with $\langle P'\rangle=0$ leads to 
\begin{equation}
    \label{eq:sigma0_rho1}
    \sigma_0=\left\langle\frac{1}{\sigma}\right\rangle^{-1}.
\end{equation}
Secondly, integrating $P'$ while using $\langle P\rangle=0$ allows to get the expression of $P$
\begin{equation}
    \label{eq:P_rho1}
    P(y)=\sigma_0\int_0^y\frac{1}{\sigma(z)}\mathrm{d}z-y-\sigma_0\int_0^1\left(\int_0^s\frac{1}{\sigma(z)}\mathrm{d}z\right)\mathrm{d}s+\frac{1}{2}.
\end{equation}
Thirdly, we now look at the coefficients in $\mathcal{E}(\mathcal{T}_0).$ When $\gamma\equiv 1$, integration by parts of $\langle\eqref{cell_pb_M}\times P - \eqref{cell_pb_P}\times M\rangle$ leads to
\begin{equation}
    \label{eq:IPPM_rho}
    \langle \sigma M' \rangle =\langle \sigma SP'\rangle-\langle \sigma S'P\rangle+\langle P^2\rangle-\frac{1}{\sigma_0}\langle \sigma P^2\rangle.
\end{equation} 
Using \eqref{eq:simpl_sigma0}, $\langle S\rangle = 0$, and \eqref{eq:Sprim_rho1}, we can then write 
\begin{equation}
    \label{first_term_DD_Erho1}
    \langle\sigma S\rangle+\langle\sigma M'\rangle  = \langle P^2\rangle. 
\end{equation}
Moreover, integration by parts of $\langle\eqref{cell_pb_N}\times P - \eqref{cell_pb_P}\times N\rangle$ leads to
\begin{equation}
    \label{eq:IPPN_rho}
    \langle \sigma N' \rangle =-v\langle PR'\rangle,
\end{equation} 
where one can show by integration of \eqref{cell_pb_R} together with $\langle R'\rangle=0$, that $R'$ is given by 
\begin{equation}
    \label{eq:Rprim_rho1}
    R'=\frac{v}{\sigma}\left(-P+\left\langle\frac{P}{\sigma}\right\rangle\left\langle\frac{1}{\sigma}\right\rangle^{-1}\right).
\end{equation}
Noting that $\frac{1}{\sigma}=\frac{P'+1}{\sigma_0}$ gives $\langle\frac{P}{\sigma}\rangle=0.$
Consequently, one gets
\begin{equation}
    \langle\sigma N'\rangle =v^2\left\langle\frac{P^2}{\sigma}\right\rangle.
\end{equation}
Finally, the effective equation \eqref{eq:def_eff_0} reads in terms of $\sigma$:
\begin{equation}
\label{eq:final_eff_eq_rho1_non_dim}
    \begin{aligned}    
     \partial_t \mathcal{T}^{(2)}-\left(\sigma_0+\eta^2v^2\left\langle\frac{P^2}{\sigma}\right\rangle\right)\partial^2_{xx}\mathcal{T}^{(2)} -2 \eta^2 v\langle P^2\rangle\partial^3_{xxx}\mathcal{T}^{(2)}-\eta^2\langle P^2\rangle\partial^3_{txx}\mathcal{T}^{(2)}=0,
     \end{aligned}
\end{equation}
with $\sigma_0$ and $P$ given by \eqref{eq:sigma0_rho1} and \eqref{eq:P_rho1}. Using the definition of the non-dimensionalized coordinates \eqref{eq:nondimparam}, we can write its dimensionalized counterpart: 
\begin{equation}
\label{eq:final_eff_eq_rho1}
    \begin{aligned}    
     \bar{\gamma}\partial_\tau\Theta^{(2)}-\left(\bar{\sigma}_0  +h^2 v_m^2  \bar{\gamma}^2 \beta_2\right)\partial^2_{XX}\Theta^{(2)} -2 h^2 v_m  \bar{\gamma}\beta_1\partial^3_{XXX}\Theta^{(2)}-h^2 \bar{\gamma}\beta_1\partial^3_{\tau XX}\Theta^{(2)}=0
     \end{aligned}
\end{equation}
with 
\begin{equation}
    \label{def_betas}
     \beta_1 = \langle P^2\rangle>0 \text{ and } \beta_2 =  \frac{1}{\sigma^\star}\left\langle\frac{P^2}{\sigma}\right\rangle>0.
\end{equation}
Integration by parts yields the energy balance 
\begin{equation}
    \label{energy_balance}
    \frac{\mathrm{d}\mathscr{E}}{\mathrm{d}\tau}=-\left(\bar{\sigma}_0  +h^2v_m^2 \bar{\gamma}^2\beta_2\right)\int_{\mathbb{R}} (\partial_X \Theta^{(2)})^2\mathrm{d}X
\end{equation}
with 
\begin{equation}
    \label{def_energy}
    \mathscr{E} =\frac{\bar{\gamma}}{2}\int_{\mathbb{R}}\left((\Theta^{(2)})^2+h^2\beta_1(\partial_X \Theta^{(2)})^2\right)\mathrm{d}X
\end{equation}
which describes a decreasing energy. One notices that except for the correction $h^2v_m^2\bar{\gamma}^2\beta_2$, the second-order terms are purely dispersive terms, which do not dissipate energy. \\
One also notices that the non-reciprocal term $\partial^3_{XXX}\Theta^{(2)}$ can be replaced by a $\partial^2_{X\tau}\Theta^{(2)}$ term up to a factor $\mathcal{O}(h)$. Therefore, one difference with the case of both parameters modulated is that the non-reciprocity occurs in the transient regime only, contrary to the convective term of \eqref{homog_order0_adim} which remains in quasi-steady heat transfer. \\
For the bilayered medium described in Section \ref{Sec:BF}, the coefficients $\beta_1$ and $\beta_2$ can be computed analytically 
\begin{equation}
    \label{beta1_2}
    \beta_1 = \frac{(\bar{\sigma}_A-\bar{\sigma}_B)^2(-1+\varphi^2)\varphi^2}{12(\bar{\sigma}_A-\bar{\sigma}_A\varphi+\bar{\sigma}_B\varphi)^2} \text{ and } \beta_2 = -\frac{(\bar{\sigma}_A-\bar{\sigma}_B)^2(-1+\varphi^2)\varphi^2}{12\bar{\sigma}_A\bar{\sigma}_B(\bar{\sigma}_A(-1+\varphi)-\bar{\sigma}_B\varphi)}, 
 \end{equation}
 and hence all the effective parameters in \eqref{eq:final_eff_eq_rho1} can be computed. \\
 Snapshots for the second-order homogenized field are presented in Figure \ref{fig:Homog2_rho1_field} for a bilayer defined by  
\begin{equation}
    \label{eq:bilayer_param_num}
    \left\lbrace
    \begin{aligned}
        &(\bar{\sigma}_A,\bar{\gamma}_A)=(500\,\mbox{ kg m s}^{-3}\mbox{ K}^{-1},\,2\times 10^6\mbox{ kg m}^{-1}\,\mbox{s}^{-2}\mbox{ K}^{-1}),\\
        &(\bar{\sigma}_B,\bar{\gamma}_B)=(1\times 10^5\mbox{ kg,\,m\,s}^{-3}\mbox{ K}^{-1},\,2\times 10^6\mbox{ kg m}^{-1}\mbox{s}^{-2}\mbox{ K}^{-1}),
    \end{aligned}
    \right.
\end{equation}
a modulation velocity $v_m=5\times 10^{-2}\mbox{ m s}^{-1}$, a volume fraction $\varphi=0.2$ and a periodicity $h=10^{-1}$m. The initial condition \eqref{GaussianInit} is chosen to be centered at $X_0=3$m, with various standard deviation $\nu$. The typical wavelength for this initial data is assumed to be about $6\nu$, so that the associated so-called small parameter of the long-wavelength homogenization regime defined in \eqref{eq:nondimparam} is $\eta = 2\pi h/(6\nu)=\pi/(30\nu)$.\\
One can see that for very small $\eta$, there is no apparent dispersive effect, which is expected since the dispersion is described only from the second order. However as we increase $\eta$ by considering smaller wavelengths, oscillations and non-reciprocity are then observed. Even if the associated values of $\eta$ are no longer negligible compared to 1, we will see on dispersion diagrams in the next section that these phenomena happen for values at which we still have a good agreement between the second-order homogenized model and the exact dispersion relation in the microstructure.

\subsection{Effective dispersion relation for $\gamma$ constant}
\begin{figure}[h]
    \subfloat[Imaginary part \label{fig:Homog2_rho1_imag}]{%
    \includegraphics[width=0.48\textwidth]{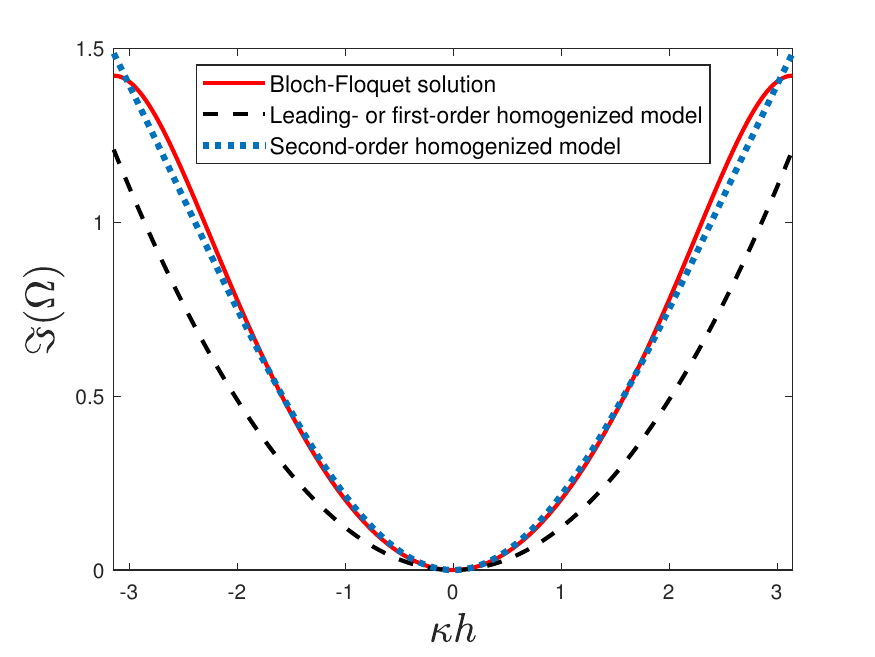}
    }
    \hfill
    \subfloat[Real part \label{fig:Homog2_rho1_real}]{%
    \includegraphics[width=0.48\textwidth]{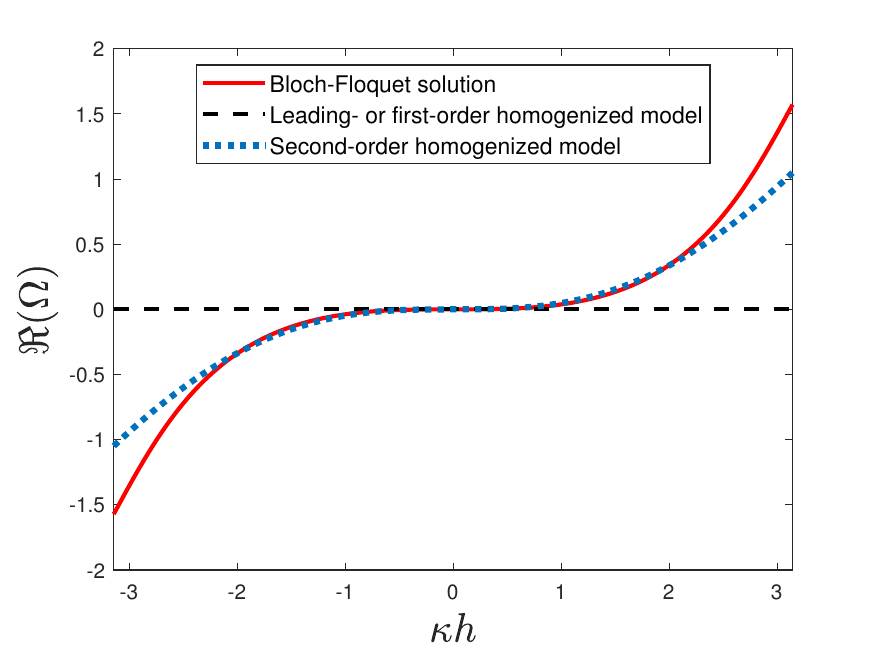}
    }\\
    \centering
    \subfloat[Complex plane\label{fig:Homog2_rho1_compl}]{%
    \includegraphics[width=0.48\textwidth]{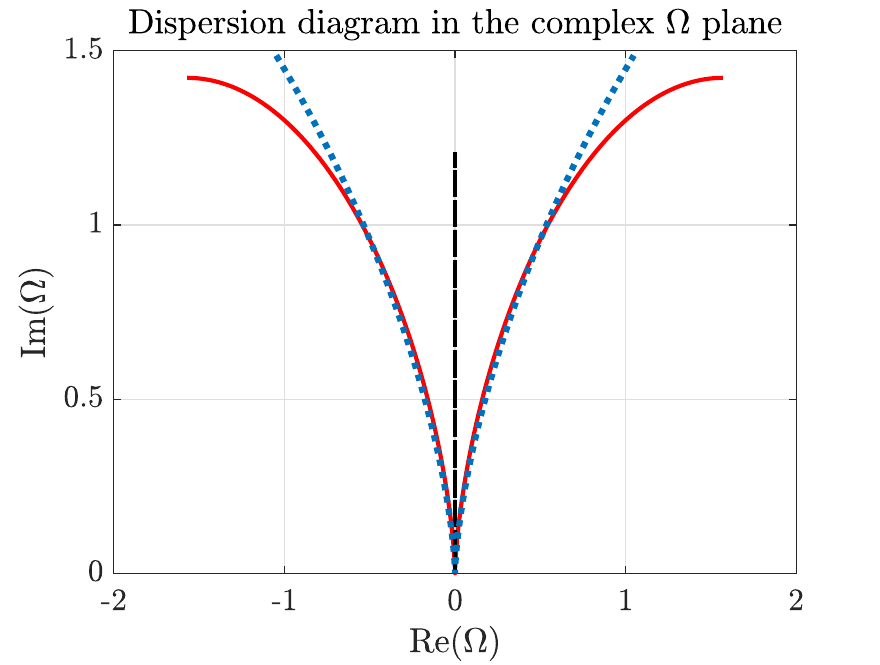}
    }
    \caption{Dispersion diagrams when only $\sigma$ is modulated (in this case Model 1 and Model 2 are the same). Exact one in plain lines, leading-order homogenized model in dashed lines, second-order homogenized model in dotted lines.}
    \label{fig:Homog2_rho1}
\end{figure}

\begin{figure}[h]
    \subfloat[Imaginary part \label{fig:Higher_branches_imag}]{%
    \includegraphics[width=0.48\textwidth]{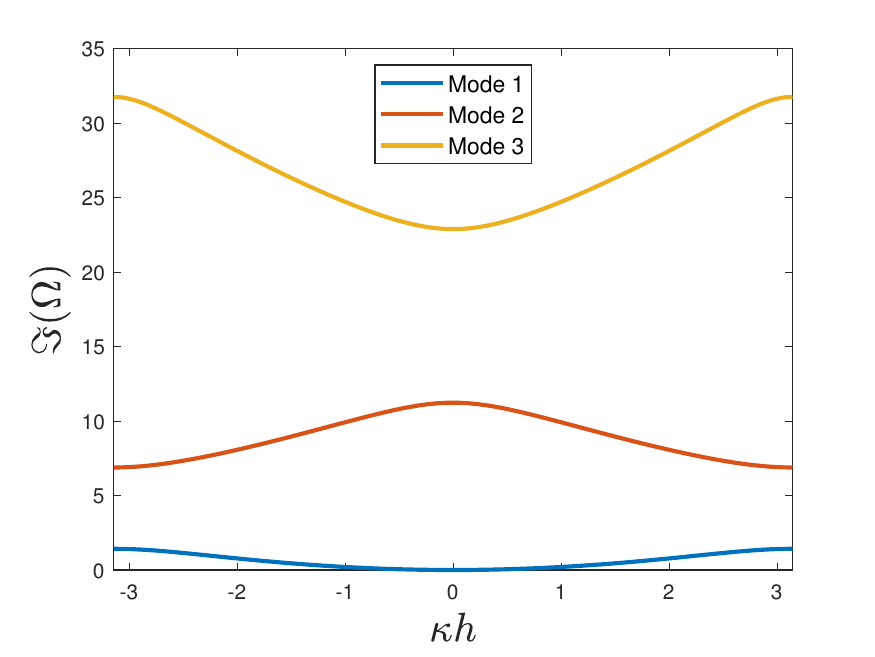}
    }
    \hfill
    \subfloat[Real part \label{fig:Higher_branches_real}]{%
    \includegraphics[width=0.48\textwidth]{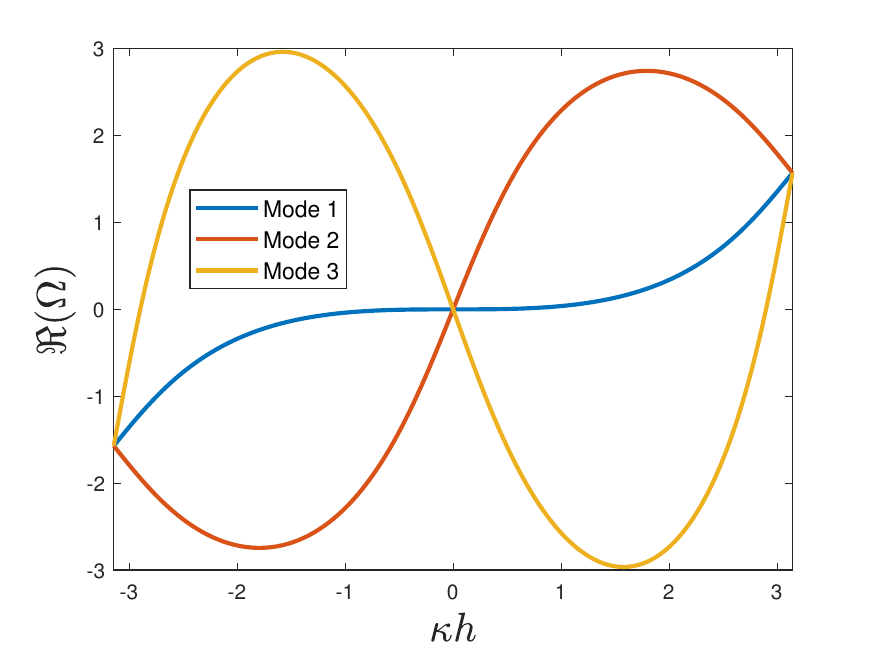}
    }\\
    \centering
    \subfloat[Complex plane\label{fig:Higher_branches_compl}]{%
    \includegraphics[width=0.48\textwidth]{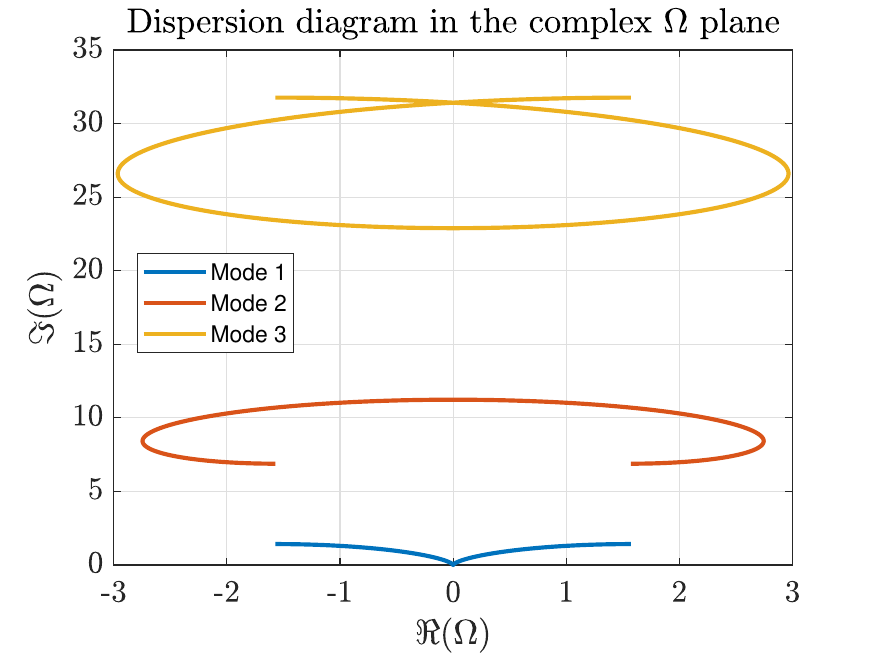}
    }
    \caption{Dispersion diagrams when only $\sigma$ is modulated  (in this case Model 1 and Model 2 are the same): first three branches obtained by Floquet-Bloch analysis.}
    \label{fig:Higher_branches}
\end{figure}

We consider a harmonic solution $\Theta^{(2)}=\mathrm{exp}(\mathrm{i}(\Omega \tau - \kappa X))$ inserted in \eqref{eq:final_eff_eq_rho1} to get the following dispersion relation: 
\begin{equation}
    \mathrm{i} \bar{\gamma}\Omega+\left(\bar{\sigma}_0  +h^2 v_m^2 \bar{\gamma}^2\beta_2\right)\kappa^2-2 h^2 v_m \bar{\gamma}\beta_1\mathrm{i}\kappa^3+h^2 \bar{\gamma}\beta_1\mathrm{i}\Omega \kappa^2=0,
\end{equation}
leading to 
\begin{equation}
\label{eq:DD_order2}
    \left\lbrace
    \begin{aligned}
        & \Re (\Omega) = \frac{2  v_m h^2\beta_1 \kappa^3}{1+h^2\beta_1  \kappa^2} \neq 0\\
        & \Im (\Omega) = \frac{\bar{\sigma}_0  +h^2 v_m^2 \bar{\gamma}^2\beta_2}{ \bar{\gamma}+h^2 \bar{\gamma}\beta_1 \kappa^2}\kappa^2
    \end{aligned}
    \right.
\end{equation}
with $\beta_1$ and $\beta_2$ defined in \eqref{def_betas}. At the leading order and first order, \eqref{eq:DD_order2} simplifies to $\Re (\Omega)=0$ and $\Im (\Omega)=\bar{\sigma}_0  \kappa^2$. 
One can see that $\Re (\Omega)$ in \eqref{eq:DD_order2} is an odd function of $\kappa$ with  $\Re (\Omega)\neq 0 $ if $\kappa\neq 0$ and $v_m\neq 0$. Moreover, $\Re'(\Omega)$ has the same sign as $v_m$. This proves that as soon as there is a modulation of the conductivity, there is propagation in the direction of the modulation. The propagation feature in the effective equation (\ref{eq:DD_order2}) is reminiscent of thermal problems with wave propagation \cite{joseph1989heat,farhat2019scattering,cassier2022active}. \\
For a bilayer, the dispersion relation is therefore given analytically by \eqref{eq:DD_order2} together with \eqref{beta1_2}. It is compared to a Taylor expansion up to the second order as $h$ tends to 0 of the dispersion relation obtained by the Floquet-Bloch analysis in Section \ref{Sec:BF}. Both expressions are exactly the same validating, therefore our second-order model.\\
The dispersion diagrams are plotted in Figure \ref{fig:Homog2_rho1}  for the same parameters as in Figure \ref{fig:Homog2_rho1_field}. The plain lines stand for the exact dispersion diagram obtained by Floquet-Bloch analysis. The black dashed lines represent the dispersion relation associated to either the leading- or first-order homogenized model obtained in the literature for the case of a sinusoidal modulation and extended in the present paper to the case of any modulation. The blue dotted line is used for the second-order homogenized model developed for the first time in this paper. One can see the better agreement of this last model as we go further from the origin, and the fact that, since $\Re(\Omega)$ is a non-zero and non-symmetric function of $\kappa$, it describes propagation of the temperature field with a non-reciprocal behaviour contrary to the lower-order ones.  \\
These dispersion diagrams also confirm the results of the simulations presented in Figure \ref{fig:Homog2_rho1_field}: very close to the origin, the real part is very small explaning the fact that the non-symmetric dispersive effects where hardly noticeable for the smallest chosen value of $\eta$. However, as one goes further away from the origin, real and imaginary parts are of the same order of magnitude, explaining stronger dispersive effects for higher values of $\eta$. Moreover, one can see that for the values of $\eta$ chosen in Figure \ref{fig:Homog2_rho1_field}, our model still presents a good fit with the exact dispersion relation; the simulations should therefore describe accurately the behaviour of the mean field in the microstructured medium. \\
The next higher-frequency solutions of the exact dispersion relation are also plotted in Figure \ref{fig:Higher_branches}. One can see from the sign of the real part that a propagation of the temperature field against the direction of the modulation is possible, as noticed in \cite{Xu2022}. However, this phenomenon occurs for much higher frequencies, meaning that no low-frequency effective model would be able to describe them. It could motivate the development of high-frequency homogenization methods for the diffusion equation \cite{Craster2010,Touboul2023}.

\section{Homogenization for Model 2}\label{Sec:Homog_advection}
It is stated in \cite{Li2022} that reciprocity is recovered when considering Model 2. However, it is obviously a limitation of leading-order homogenization since, when only $\bar{\sigma}$ is modulated, the models \eqref{dim_1D_eq} and \eqref{dim_1D_eq_adv} are the same and we have proved in Section \ref{Sec:Homog} that high-order homogenization allows to describe a non-reciprocal behaviour. \\
The same approach as the one followed for Model 1 can be conducted on \eqref{dim_1D_eq_adv}, and for the sake of completeness the main results are summarized in this section. 
The non-dimensionalized counterpart of \eqref{dim_1D_eq_adv} is: 
\begin{equation}
    \label{non_dim_1D_eq_adv}
      \partial_x\left( {\sigma\left(\frac{x-vt}{\eta}\right)}\partial_x {T_\eta} \right)+{f} = c\rho\left(\frac{x-vt}{\eta}\right) \partial_t {T_\eta}+c\left[ {\rho}\left(\frac{x-vt}{\eta}\right)-{\rho}_0\right]v\partial_x {T_\eta}
\end{equation}
with $\rho_0=\langle \rho\rangle$ and together with continuity conditions for $T_\eta$ and $\sigma \partial_x T_\eta$. \\
The leading-order field is still a macroscopic one $T_0(x,y,t)=\mathcal{T}_0(x,t)$ given by the following effective equation: 
\begin{equation}
    \label{eq:eff_order0_adv}
    \gamma_0 \partial_t \mathcal{T}_0 = \sigma_0 \partial^2_{xx}\mathcal{T}_0+f
\end{equation}
with
\begin{equation}
    \label{eff_param_order0_adv}
    \gamma_0=c\langle \rho\rangle \text{ and } \sigma_0=\langle\sigma(P'+1)\rangle
\end{equation}
as in \eqref{effective_parameters_zero_order}, with $P$ still given by \eqref{cell_pb_P}. This effective equation is indeed reciprocal as stated in \cite{Li2022}. \\
The first-order field can be written as 
\begin{equation}
    \label{adv_field_T1}
    T_1(x,y,t)=\mathcal{T}_1(x,t)+P(y)\partial_x\mathcal{T}_0(x,t)
\end{equation}
where the first-order macroscopic field satisfies 
\begin{equation}
    \label{eq:eff_order1_adv}
    \gamma_0 \partial_t \mathcal{T}_1= \sigma_0 \partial^2_{xx}\mathcal{T}_1
\end{equation}
which is again reciprocal. \\
Finally, the second-order field reads
\begin{equation}
    \label{eq:second_order_decompo_adv}
    T_2(x,y,t)=\mathcal{T}_2(x,t)+P(y)\partial_x\mathcal{T}_1(x,t)+R_\mathrm{adv}(y)\partial_x\mathcal{T}_0(x,t)+S(y)\partial_{t}\mathcal{T}_0(x,t)+A(y)f(x,t)
\end{equation}
where $S$, and $A$ satisfy the same cell problems as for Model 1, i.e. \eqref{cell_pb_S} and \eqref{cell_pb_A}, respectively, while $R_\mathrm{adv}$ is solution of a new cell problem given by \eqref{cell_pb_R_adv}. 

The second-order mean field $\mathcal{T}_2$ satisfies the following effective equation: 
\begin{equation}
    \label{eff_order2_adv}
    \begin{aligned}
    \gamma_0\partial_t\mathcal{T}_2-\sigma_0\partial^2_{xx}\mathcal{T}_2&=\left[-\frac{\sigma_0}{\rho_0}\langle\rho S\rangle +\langle\sigma S\rangle +\langle\sigma M'\rangle \right]\partial^3_{txx}\mathcal{T}_0+\left[-cv\langle\rho R_\mathrm{adv}\rangle +\langle\sigma N_\mathrm{adv}'\rangle  \right]\partial^2_{xx}\mathcal{T}_0\\
    &+ \left[-\frac{\sigma_0}{\rho_0}\langle\rho R_\mathrm{adv}\rangle -v\frac{\sigma_0}{\rho_0}\langle \rho S\rangle +\langle\sigma R_\mathrm{adv}\rangle+\langle\sigma L_\mathrm{adv}'\rangle \right]\partial^3_{xxx}\mathcal{T}_0 \\ 
    & +\left[-\frac{1}{\rho_0}\langle\rho S\rangle-c\langle\rho A\rangle \right]\partial_t f +\left[\langle\sigma C'\rangle+\langle \sigma A\rangle \right]\partial^2_{xx}f \\
    & +\left[-\frac{1}{\rho_0}\langle\rho R_\mathrm{adv}\rangle-cv\langle\rho A\rangle+\langle \sigma B_\mathrm{adv}'\rangle -\frac{v}{\rho_0}\langle\rho S\rangle\right]\partial_{x}f
     \end{aligned}
\end{equation}
where $M$ and $C$ satisfy \eqref{cell_pb_M} and \eqref{cell_pb_C}, respectively, as for Model 1, while $L_\mathrm{adv}$, $N_\mathrm{adv}$, and $B_\mathrm{adv}$ satisfy the cell problems given in  \ref{app:cell_Model2}. 

The non-reciprocal behaviour is therefore encapsulated in the coefficient 
\begin{equation}
    \label{eq:Ncoeff_adv}
    \begin{aligned}
    \mathcal{N}_\mathrm{adv}=-\frac{\sigma_0}{\rho_0}\langle\rho R_\mathrm{adv}\rangle -v\frac{\sigma_0}{\rho_0}\langle \rho S\rangle +\langle\sigma R_\mathrm{adv}\rangle+\langle\sigma L_\mathrm{adv}'\rangle.
    \end{aligned}
\end{equation}
Assuming that $\rho$ can be written as 
\begin{equation}
    \label{decompo_rho}
    \rho(y) = \rho_0+\delta\rho(y) 
\end{equation}
and chosing $\sigma^\star$ such that 
$$\left\langle\frac{1}{\sigma}\right\rangle=1$$
one can prove that $ \mathcal{N}_\mathrm{adv}  $ reduces to 
\begin{equation}
    \label{eq:Ncoeff_adv_final}
    \begin{aligned}
    \mathcal{N}_\mathrm{adv}=\frac{2c\sigma_0v}{\rho_0}\left[ \left\langle\left(\int_0^y\delta\rho(z)\mathrm{d}z\right)^2\frac{1}{\sigma}\right\rangle-\left\langle\frac{1}{\sigma}\int_0^y\delta\rho(z)\mathrm{d}z\right\rangle^2\right]+2cv\langle\rho P^2\rangle-4cv\left\langle P\int_0^y \delta \rho(z)\mathrm{d}z\right\rangle.
    \end{aligned}
\end{equation}
We do not prove that this term is always non-zero as we did for Model 1. However, in all the cases investigated numerically it was observed to be so. We also give in the next three sections examples of families of choice parameters for which $ \mathcal{N}_\mathrm{adv}$ can be proved to be non-zero. In all these cases, propagation and non-reciprocity is therefore described by the second-order model. 
\subsection{Only $\sigma$ modulated}
As said previously, if only $\sigma$ is modulated Models 1 and 2 are the same, and non-reciprocity always holds at the second order. Indeed, in that case, $\delta\rho=0$ and $\mathcal{N}_\mathrm{adv}$ can be written as 
\begin{equation}
    \label{eq:Ncoeff_adv_final_rho1}
    \begin{aligned}
    \mathcal{N}_\mathrm{adv}=2cv\rho_0\langle P^2\rangle
    \end{aligned}
\end{equation}
which is non-zero since $P\neq 0$ as soon as $\sigma$ is non-constant. This case has therefore been illustrated in the previous section. 
\subsection{Only $\rho$ modulated}
\begin{figure}[htbp]
    \subfloat[Imaginary part \label{fig:Homog2_rho1_imag}]{%
    \includegraphics[width=0.48\textwidth]{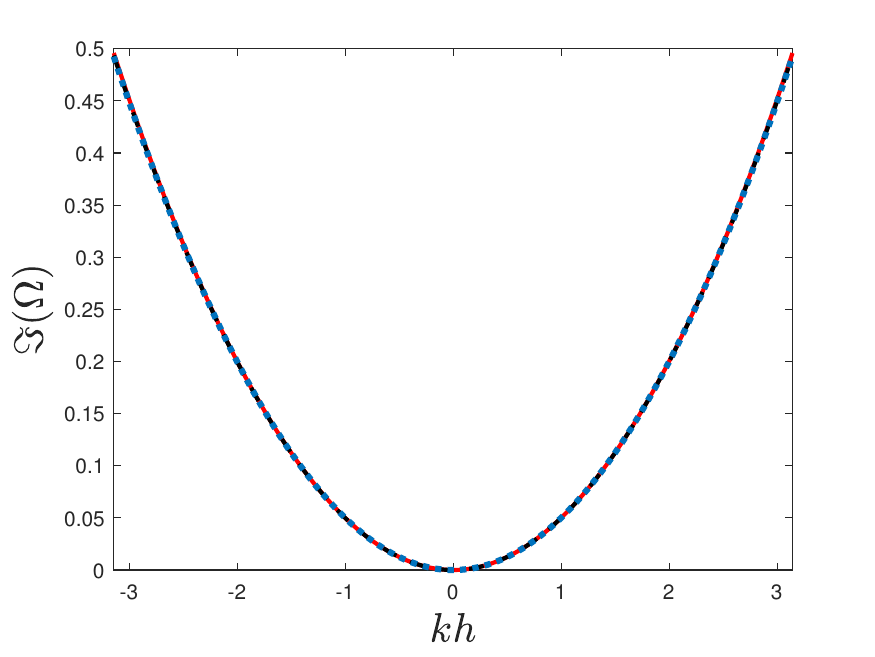}
    }
    \hfill
    \subfloat[Real part \label{fig:Homog2_rho1_real}]{%
    \includegraphics[width=0.48\textwidth]{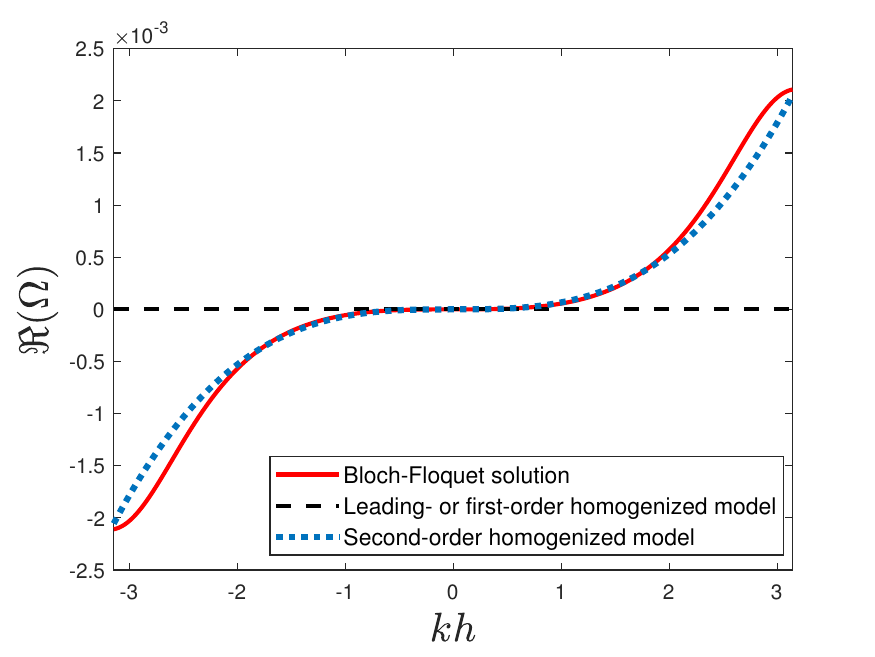}
    }\\
    \centering
    \subfloat[Complex plane\label{fig:Homog2_rho1_compl}]{%
    \includegraphics[width=0.48\textwidth]{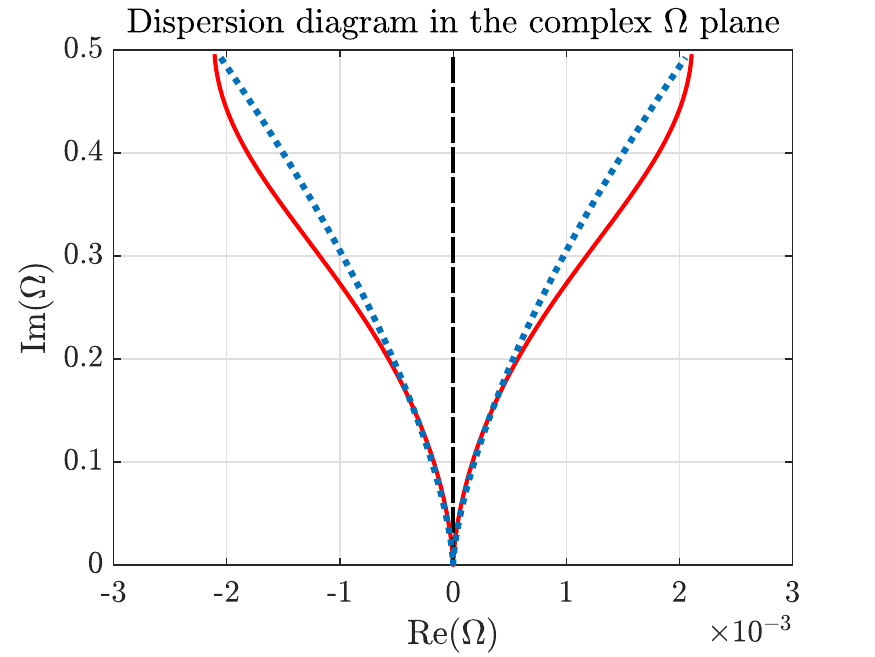}
    }
    \caption{Dispersion diagrams when only $\rho$ is modulated (Model 2). Exact one in plain lines, leading-order homogenized model in dashed lines, second-order homogenized model in dotted lines.}
    \label{fig:Homog2_rho1_adv}
\end{figure}

When $\sigma=1$, $P=0$ and $ \mathcal{N}_\mathrm{adv}  $ simplifies to 
\begin{equation}
    \label{eq:Ncoeff_adv_final_sigma1}
    \begin{aligned}
    \mathcal{N}_\mathrm{adv}=\frac{2c v}{\rho_0}\left[ \left\langle\left(\int_0^y\delta\rho(z)\mathrm{d}z\right)^2\right\rangle-\left\langle\int_0^y\delta\rho(z)\mathrm{d}z\right\rangle^2\right].
    \end{aligned}
\end{equation}
The equality $\mathcal{N}_\mathrm{adv}=0$ in \eqref{eq:Ncoeff_adv_final_sigma1} would mean that $\int_0^y\delta\rho(z)\mathrm{d}z$ is constant (case of equality in Cauchy-Schwarz inequality) and therefore $\delta\rho=0$, which is incompatible with the fact that $\rho$ is modulated. Therefore, here again, we know that the term responsible for non-reciprocity and propagation is always non-zero. \\
In this case, the effective equation reads in dimensionalized coordinates for $F=0$:
\begin{equation}
\label{eq:final_eff_eq_adv}
    \begin{aligned}    
     \bar{\gamma}_0\partial_\tau\Theta^{(2)}-\left(\bar{\sigma}  +h^2 v_m^2 \frac{\bar{c}^2}{\bar{\sigma}}\beta_3\right)\partial^2_{XX}\Theta^{(2)} -2 h^2 v_m  \frac{\bar{c}}{\bar{\rho}_0}\beta_3\partial^3_{XXX}\Theta^{(2)}-h^2 \frac{\bar{c}}{\bar{\rho}_0}\beta_3\partial^3_{\tau XX}\Theta^{(2)}=0
     \end{aligned}
\end{equation}
where $\beta_3$ is defined by
\begin{equation}
    \label{eq:beta3}
    \beta_3 = (\rho^\star)^2 \langle (S')^2\rangle 
\end{equation}
and in the case of a bilayered laminate can be written as
\begin{equation}
    \label{eq:beta3_bilayer}
    \beta_3 = \frac{1}{12}(\bar{\rho}_A-\bar{\rho}_B)^2(-1+\phi)^2\phi^2.
\end{equation}
\begin{equation}
    \mathrm{i} \bar{\gamma}_0\Omega+\left(\bar{\sigma}  +h^2 v_m^2 \frac{\bar{c}^2}{\bar{\sigma}}\beta_3\right)\kappa^2-2 h^2 v_m  \frac{\bar{c}}{\bar{\rho}_0}\beta_3\mathrm{i}\kappa^3+h^2 \frac{\bar{c}}{\bar{\rho}_0}\beta_3\mathrm{i}\Omega \kappa^2=0.
\end{equation}
This expression is validated with respect to a formal Taylor expansion in $h$ of the exact dispersion relation obtained with Bloch-Floquet analysis. A comparison of the dispersion diagrams is presented in Figure \ref{fig:Homog2_rho1_adv} for a bilayer defined by $c=1000\mbox{ m}^2\mbox{ K}^{-1}\mbox{ s}^{-2}$ together with
\begin{equation}
    \label{eq:bilayer_param_num_adv}
    \left\lbrace
    \begin{aligned}
        &(\bar{\sigma}_A,\bar{\rho}_A)=(190\mbox{ kg m s}^{-3}\mbox{ K}^{-1},\,2\times 10^6\mbox{ kg m}^{-1}\,\mbox{ s}^{-2}\mbox{ K}^{-1}),\\
        &(\bar{\sigma}_B,\bar{\rho}_B)=(190\mbox{ kg m s}^{-3}\mbox{ K}^{-1},\,1.4\times 10^6\mbox{ kg m}^{-1}\,\mbox{ s}^{-2}\mbox{ K}^{-1}),
    \end{aligned}
    \right.
\end{equation}
a modulation velocity $v_m=5\times 10^{-3}\mbox{ m s}^{-1}$, a volume fraction $\varphi=0.2$ and a periodicity $h=5\times 10^{-2}\mbox{ m}$. One can observe a good description of the propagative part missed by the leading order.
\subsection{Bilayered material}
For a bilayered material, whatever the choice of the parameters and the modulation are, $\mathcal{N}_\mathrm{adv}$ can be expressed as a non-zero polynominal of degree 4 in the volume fraction $\varphi$. Consequently, it is zero in at most 4 values on $(0,1)$. Consequently, the term almost never vanishes and there are infinitely many configurations for which non-reciprocity and propagation occur. 



\section{Conclusion}
The articles \cite{Torrent2018} and \cite{Xu2022} derived an effective equation from a Floquet-Bloch expansion in the time-modulated diffusion at the leading and first orders respectively; here we push this analysis to the second order. We uncover non-reciprocal features in the effective equation at the second order when a single parameter is modulated, in stark contrast with \cite{Torrent2018,Xu2022} for which non-reciprocity requires time-modulation of both parameters. Our approach based on two-scale asymptotic expansions \cite{Bensoussan,Bakhvalov1989,cioranescu1999,Conca1997} is more general and allows to consider arbitrary laminates with source terms. Our theoretical and numerical results shed a new light on diffusion processes in time-modulated media, and complement our former work on second-order homogenization of the time-modulated wave equation \cite{Touboul2024}. We hope that our results will foster experimental efforts in non-reciprocal effects notably in mass, particle and light diffusion processes for which only one parameter can be time-modulated. We finally note that the marked differences between second-order homogenization of the time-modulated wave equation \cite{Touboul2024} and the present work, suggest the theory of field patterns introduced for the wave equation in \cite{milton2017field}, would be worth investigating for the diffusion equation. \\

\noindent {\bf Acknowledgements} SG and RVC were funded by UK Research and Innovation (UKRI) under the UK government’s Horizon Europe funding guarantee [grant number 10033143].
RVC thanks UKRI and EPSRC for support through grant number EP/Y015673/1.

\appendix 
\section{Cell problems}
\subsection{Leading order for Model 1}\label{app:cell_leading}
The correctors $P$ and $Q$ used to write $T_1$ in \eqref{form_T1} are defined by  
\begin{equation}
    \label{cell_pb_P}
    \left\lbrace
    \begin{aligned}
    &\partial_y \left[\sigma(P'+1)\right]=0 \\
    & \text{ continuity of $P$ and $\sigma(P'+1)$} \\
    & \text{1-periodicity of } P \text{ and } \langle P \rangle = 0 ,
    \end{aligned}
    \right.
\end{equation}
and
\begin{equation}
    \label{cell_pb_Q}
    \left\lbrace
    \begin{aligned}
    &\partial_y \left[\sigma Q'+v\gamma\right]=0 \\
    &\text{continuity of $Q$ and $\sigma Q'+v\gamma$} \\
    &  \text{1-periodicity of } Q \text{ and } \langle Q \rangle = 0, 
    \end{aligned}
    \right.
\end{equation}
so that \eqref{eq:identif_etam1} and the continuity conditions for $T_1$ and $-a_0+v\gamma T_0$ are satisfied. We notice for the next orders that
\begin{equation}
\label{eq:simpl_sigma0}
   \sigma(P'+1)=\sigma_0 
\end{equation}
with $\sigma_0$ the effective parameter defined in \eqref{effective_parameters_zero_order}. Furthermore, integration by parts of $\langle\eqref{cell_pb_Q}\times P - \eqref{cell_pb_P}\times Q\rangle$ leads to
\begin{equation}
    \label{eq:IPP1}
    \langle\sigma Q'\rangle = v\langle\gamma P'\rangle
\end{equation}
\subsection{First order for Model 1}\label{app:cell_first}
The correctors $R$, $S$, $V$ and $A$ allowing to write $T_2$ in \eqref{form_T2} are given by the following cell problems
\begin{equation}
    \label{cell_pb_R}
    \left\lbrace
    \begin{aligned}
    &\partial_y \left[\sigma R'+v\gamma P+\sigma Q -\frac{W_0}{\sigma_0}\sigma P\right]=-\sigma Q'+W_0 \\
    &  \text{continuity of $R$ and $\sigma R'+v\gamma P+\sigma Q -\frac{W_0}{\sigma_0}\sigma P$} \\
    & \text{1-periodicity of } R \text{ and } \langle R \rangle = 0, 
    \end{aligned}
    \right.
\end{equation}
\begin{equation}
    \label{cell_pb_S}
    \left\lbrace
    \begin{aligned}
    &\partial_y \left[\sigma S'+\frac{\gamma_0}{\sigma_0}\sigma P\right]=\gamma-\gamma_0 \\
    &  \text{continuity of $S$ and $\sigma S'+\frac{\gamma_0}{\sigma_0}\sigma P$} \\
    & \text{1-periodicity of } S \text{ and } \langle S \rangle = 0, 
    \end{aligned}
    \right.
\end{equation}
\begin{equation}
    \label{cell_pb_V}
    \left\lbrace
    \begin{aligned}
    &\partial_y \left[\sigma V'+v\gamma Q\right]=0 \\
    &  \text{continuity of $V$ and $\sigma V'+v\gamma Q$} \\
    & \text{1-periodicity of } V \text{ and } \langle V \rangle = 0, 
    \end{aligned}
    \right.
\end{equation}
\begin{equation}
    \label{cell_pb_A}
    \left\lbrace
    \begin{aligned}
    &\partial_y \left[\sigma A'-\frac{1}{\sigma_0}\sigma P\right]=0\\
    &  \text{continuity of $A$ and $\sigma A'-\frac{1}{\sigma_0}\sigma P$} \\
    & \text{1-periodicity of } A \text{ and } \langle A \rangle = 0, 
    \end{aligned}
    \right.
\end{equation}
so that \eqref{syst:T_2} and continuity for $T_2$ and $-a_1+v\gamma T_1$ is satisfied. \\
Integration by parts of $\langle\eqref{cell_pb_R}\times P - \eqref{cell_pb_P}\times R\rangle$ leads to
\begin{equation}
    \label{eq:IPP2}
    \langle \sigma R' \rangle =v\langle\gamma PP'\rangle+\langle\sigma QP'\rangle+\frac{W_0}{\sigma_0}\langle\sigma P\rangle-\langle\sigma Q'P\rangle.
\end{equation} 
Integration by parts of $\langle\eqref{cell_pb_S}\times P - \eqref{cell_pb_P}\times S\rangle$ leads to
\begin{equation}
    \label{eq:IPP4}
    \langle \sigma S' \rangle =\langle\gamma P\rangle -\frac{\gamma_0}{\sigma_0}\langle \sigma P\rangle .
\end{equation}
Integration by parts of $\langle\eqref{cell_pb_V}\times P - \eqref{cell_pb_P}\times V\rangle$ leads to
\begin{equation}
    \label{eq:IPP6}
    \langle \sigma V' \rangle =v\langle\gamma QP'\rangle.
\end{equation} 
Integration by parts of $\langle\eqref{cell_pb_A}\times P - \eqref{cell_pb_P}\times A\rangle$ leads to
\begin{equation}
    \label{eq:IPP8}
    \langle \sigma A' \rangle =\frac{1}{\sigma_0}\langle \sigma P\rangle .
\end{equation}
\subsection{Second order for Model 1} \label{app:cell_second}
The correctors $L$, $M$, $N$, $O$, $B$ and $C$ allowing to write $T_3$ in \eqref{form_T3} are given by the following cell problems
\begin{equation}
    \label{cell_pb_L}
    \left\lbrace
    \begin{aligned}
    &\partial_y \left[\sigma L'+\sigma R+\frac{\sigma_0}{\gamma_0}v\gamma S+\mathcal{C}_1\sigma P\right]=-\sigma(Q+R')+\frac{W_0}{\sigma_0}\sigma P+\frac{\sigma_0}{\gamma_0}\gamma Q-\sigma_0\mathcal{C}_1\\
    &  \text{continuity of $L$ and $\sigma L'+\sigma R+\frac{\sigma_0}{\gamma_0}v\gamma S+\mathcal{C}_1\sigma P$} \\
    & \text{1-periodicity of } L \text{ and } \langle L \rangle = 0, 
    \end{aligned}
    \right.
\end{equation}
\begin{equation}
    \label{cell_pb_N}
    \left\lbrace
    \begin{aligned}
    &\partial_y \left[\sigma (N'+V+\mathcal{C}_2 P) +v\gamma\left(R+\frac{W_0}{\gamma_0} S\right)\right]=-\sigma V'+\frac{W_0}{\gamma_0}\gamma Q-\sigma_0\mathcal{C}_2\\
    &  \text{continuity of $N$ and $\sigma (N'+V+\mathcal{C}_2 P) +v\gamma\left(R+\frac{W_0}{\sigma_0} S\right)$} \\
    & \text{1-periodicity of } N \text{ and } \langle N \rangle = 0, 
    \end{aligned}
    \right.
\end{equation}
\begin{equation}
    \label{cell_pb_M}
    \left\lbrace
    \begin{aligned}
    &\partial_y \left[\sigma M'+\sigma S\right]=-\sigma S'+\gamma P -\frac{\gamma_0}{\sigma_0}\sigma P\\
    &  \text{continuity of $M$ and $\sigma M'+\sigma S$} \\
    & \text{1-periodicity of } M \text{ and } \langle M \rangle = 0, 
    \end{aligned}
    \right.
\end{equation}
\begin{equation}
    \label{cell_pb_O}
    \left\lbrace
    \begin{aligned}
    &\partial_y \left[\sigma O'+v\gamma V\right]=0\\
    &  \text{continuity of $O$ and $\sigma O'+v\gamma V$} \\
    & \text{1-periodicity of } O \text{ and } \langle O \rangle = 0, 
    \end{aligned}
    \right.
\end{equation}
\begin{equation}
    \label{cell_pb_B}
    \left\lbrace
    \begin{aligned}
    &\partial_y \left[\sigma \left(B'+\frac{\langle\gamma Q\rangle}{\gamma_0\sigma_0}P \right)+v\gamma \left(A+\frac{1}{\gamma_0}S\right)\right]=\frac{\gamma}{\gamma_0}Q-\frac{\langle \gamma Q\rangle}{\gamma_0}\\
    &  \text{continuity of $B$ and $\sigma \left(B'+\frac{\langle\gamma Q\rangle}{\gamma_0\sigma_0}P \right)+v\gamma \left(A+\frac{1}{\gamma_0}S\right)$} \\
    & \text{1-periodicity of } B \text{ and } \langle B \rangle = 0, 
    \end{aligned}
    \right.
\end{equation}
\begin{equation}
    \label{cell_pb_C}
    \left\lbrace
    \begin{aligned}
    &\partial_y \left[\sigma C'+\sigma A\right]=-\sigma A'+\frac{1}{\sigma_0}\sigma P\\
    &  \text{continuity of $C$ and $\sigma C'+\sigma A$} \\
    & \text{1-periodicity of } C \text{ and } \langle C \rangle = 0. 
    \end{aligned}
    \right.
\end{equation}
\subsection{New cell problems for Model 2} \label{app:cell_Model2}
When performing homogenization for Model 2 in Section \ref{Sec:Homog_advection}, some correctors $R_\mathrm{adv}$, $L_\mathrm{adv}$, $N_\mathrm{adv}$ and $B_\mathrm{adv}$ are different from the ones given previously for Model 1. More precisely, they are given by the following cell problems
\begin{equation}
    \label{cell_pb_R_adv}
    \left\lbrace
    \begin{aligned}
    &\partial_y \left[\sigma R_\mathrm{adv}'+vc\rho_0 P\right]=\left[\rho-\rho_0\right]cv\\
    &  \text{continuity of $ R_\mathrm{adv}$ and $\sigma R'_\mathrm{adv}$} \\
    & \text{1-periodicity of } R_\mathrm{adv} \text{ and } \langle R_\mathrm{adv}\rangle = 0,
    \end{aligned}
    \right.
\end{equation}
\begin{equation}
    \label{cell_pb_L_adv}
    \left\lbrace
    \begin{aligned}
    &\partial_y \left[\sigma L_\mathrm{adv}'+v\sigma_0 S+\sigma R_\mathrm{adv}\right]=\left[\rho-\rho_0\right]cvP-\sigma R_\mathrm{adv}'\\
    &  \text{continuity of $ L_\mathrm{adv}$ and $\sigma L_\mathrm{adv}'+\sigma R_\mathrm{adv}$} \\
    & \text{1-periodicity of } L_\mathrm{adv} \text{ and } \langle L_\mathrm{adv}\rangle = 0, 
    \end{aligned}
    \right.
\end{equation}
\begin{equation}
    \label{cell_pb_N_adv}
    \left\lbrace
    \begin{aligned}
    &\partial_y \left[\sigma N_\mathrm{adv}'+\rho_0c v R_\mathrm{adv}\right]=0\\
    &  \text{continuity of $ N_\mathrm{adv}$ and $\sigma N_\mathrm{adv}'$} \\
    & \text{1-periodicity of } N_\mathrm{adv} \text{ and } \langle N_\mathrm{adv}\rangle = 0, 
    \end{aligned}
    \right.
\end{equation}
\begin{equation}
    \label{cell_pb_B_adv}
    \left\lbrace
    \begin{aligned}
    &\partial_y \left[\sigma B_\mathrm{adv}'\right]=-v(\rho_0 c A' +S')\\
    &  \text{continuity of $ B_\mathrm{adv}$ and $\sigma B_\mathrm{adv}'$} \\
    & \text{1-periodicity of } B_\mathrm{adv} \text{ and } \langle B_\mathrm{adv}\rangle = 0.
    \end{aligned}
    \right.
\end{equation}











\bibliographystyle{RS}

\bibliography{sample}
\end{document}